\documentclass[11pt,letterpaper]{article}
\usepackage[margin=1in]{geometry}
\usepackage{empheq}
\usepackage{amsmath,amssymb,amsthm,mathrsfs,bm}
\allowdisplaybreaks
\usepackage{siunitx}
\usepackage{orcidlink}
\usepackage[title]{appendix}

\usepackage{titling}
\thanksmarkseries{arabic}
\usepackage{authblk}

\usepackage[numbers,sort]{natbib}
\usepackage[normalem]{ulem}

\theoremstyle{definition}

\numberwithin{equation}{section}
\numberwithin{remark}{section}
\numberwithin{proposition}{section}
\numberwithin{corollary}{section}

\newcommand\redout{\bgroup\markoverwith
	{\textcolor{red}{\rule[0.5ex]{2pt}{0.8pt}}}\ULon}

\DeclareMathOperator{\sgn}{sgn}
\newcommand{\re}{\mathrm{e}}
\newcommand{\ri}{\mathrm{i}}
\newcommand{\rd}{\mathrm{d}}
\newcommand{\Bo}{\text{Bo}}

\newcommand{\me}{{\mathrm{e}}}
\newcommand{\dif}{\mathrm{d}}

\newcommand{\D}{\displaystyle}

\graphicspath{{figures/}}

\newcommand*{\email}[1]{\href{mailto:#1}{\nolinkurl{#1}} } 

\usepackage{hyperref}

\title{Sixth-order parabolic equation on an interval: Eigenfunction expansion, Green's function, and intermediate asymptotics for a finite thin film with elastic resistance}

\author[1,*]{Nectarios C.\ Papanicolaou\,\orcidlink{0000-0001-5739-1143}}
\affil[1]{\small\it Department of Computer Science, University of Nicosia, \authorcr\it 46 Makedonitissas Avenue, CY-2417 Nicosia, Cyprus}
\affil[*]{\email{papanicolaou.n@unic.ac.cy}}
\author[2,**]{Ivan C.\ Christov\,\orcidlink{0000-0001-8531-0531}}
\affil[2]{\small\it School of Mechanical Engineering, Purdue University, West Lafayette, Indiana 47907, USA}
\affil[**]{Corresponding author; \email{christov@purdue.edu}}

\date{\today}

\begin{document}

\maketitle

\begin{abstract}
A linear sixth-order partial differential equation (PDE) of ``parabolic'' type describes the dynamics of thin liquid films beneath surfaces with elastic bending resistance when deflections from the equilibrium film height are small. On a finite domain, the associated sixth-order eigenvalue problem is self-adjoint for the boundary conditions corresponding to a thin film in a closed trough, and the eigenfunctions form a complete orthonormal set. Using these eigenfunctions, we derive the Green's function for the governing sixth-order PDE on a finite interval and compare it to the known infinite-line solution. Further, we propose a Galerkin spectral method based on the constructed sixth-order eigenfunctions and their derivative expansions. The system of ordinary differential equations for the time-dependent expansion coefficients is solved by standard numerical methods. The numerical approach is applied to versions of the governing PDE with a second-order spatial derivative (in addition to the sixth-order one), which arises from gravity acting on the film. In the absence of gravity, we demonstrate the self-similar intermediate asymptotics of initially localized disturbances on the film surface, at least until the disturbances ``feel'' the finite boundaries, and show that the derived Green's function is an attractor for such solutions. In the presence of gravity, we use the proposed Galerkin numerical method to demonstrate that self-similar behavior persists, albeit for shortened intervals of time, even for large values of the gravity-to-bending ratio.\\[1mm]
\end{abstract}



\section{Introduction}\label{sec:Intro}

Sixth-order parabolic equations arise in a variety of physical contexts, from semiconductor manufacturing \cite{King1989} to geological intrusions \cite{Michaut2011,Bunger2011,Bunger2011b,Thorey2014}. In particular, the behavior of several classes of solutions to sixth-order parabolic equations has been analyzed in the context of single-phase, Newtonian thin liquid films spreading underneath a surface with elastic bending resistance \cite{Huang2002,Flitton2004,Hosoi2004,Hewitt2015,Carlson2016,Kodio2017,Pedersen2019,Peng2020,Pedersen2021,Saeter2024}, with applications to microfluidics \cite{Boyko2019,Boyko2020}, flow actuation \cite{Rubin2017, Boyko2019}, and impact mitigation \cite{Tulchinsky2016}. To avoid boundary conditions and simplify the analysis, many previous studies \cite{Hewitt2015,Tulchinsky2016,Rubin2017,Peng2020,Pedersen2021,Saeter2024} considered infinite domains, for example, using planar polar coordinates. Others used periodic boundary conditions \cite{Huang2002,Carlson2016}. On the other hand, understanding the effect of boundary conditions is of general interest \cite{Flitton2004,Hosoi2004} and particularly important for microfluidics applications \cite{Kodio2017,Boyko2019,Boyko2020,Gabay2023}. Similarly, the phenomenon of a thin film with elastic resistance on a finite domain may be familiar to many based on daily life---the surface of heated milk forms a skin, known as a \emph{lactoderm} composed of various proteins \cite{Martins2019}. 

In particular, \citet{Pedersen2021} studied the leveling of thin films with elastic interface resistance. They considered a 3D axisymmetric scenario, which leads to a one-dimensional partial differential equation (PDE) in polar coordinates. They demonstrated that a suitable Green's function captures the intermediate asymptotics of evolution from \emph{any} initial condition on the infinite domain. Notably, these intermediate asymptotics are \emph{self-similar} in a way dictated by the analytical construction of the Green's function. Importantly, \citet{Pedersen2021} noted that for ``bending-driven elastohydrodynamic flows, the literature is scarce compared with its capillary analog,'' despite the asymptotic adjustment being governed by a linear parabolic initial-boundary-value problem, which can be treated by the classical methods of applied mathematics. Indeed, we are not aware of previous work on the finite-interval problem for a thin film with an elastic interface, which motivates the present work. Specifically, we reconsider the analysis of \citet{Pedersen2021} of the self-similar intermediate asymptotics, but on a \emph{finite} interval. We use eigenfunction expansions to derive the Green's function for a two-dimensional thin film, which also leads to a one-dimensional (1D) PDE, but now in Cartesian coordinates (compared to polar coordinates in \cite{Pedersen2021}), on a finite interval.

\begin{figure}
	\centering
	\includegraphics[width=0.7\textwidth]{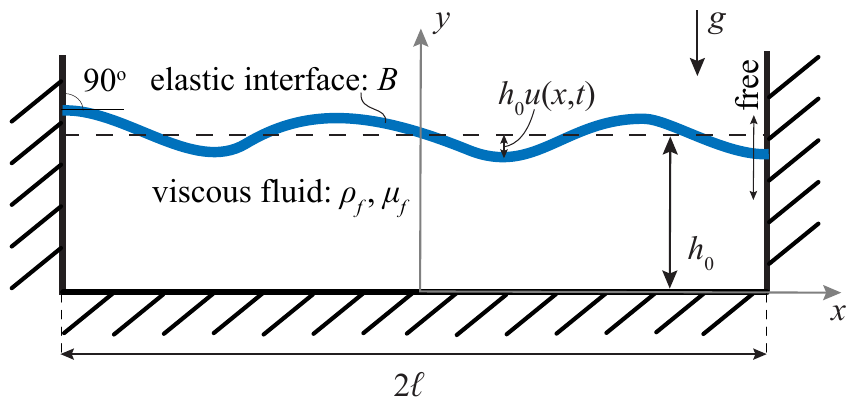}
	\caption{Schematic of the thin film with elastic resistance being considered. The interface has out-of-plane bending rigidity $B$, negligible in-plane tension, and negligible mass. The viscous fluid beneath the interface has density $\rho_f$ and dynamic viscosity $\mu_f$. The equilibrium height of the film (when the interface is flat) is $h_0$, and the dimensionless small vertical displacement of the film about this level is denoted as $u(x,t)$. The thin film is restricted to a finite domain (closed trough) of axial distance $2\ell$. The gravitational acceleration is $g$ in the $-y$ direction. Figure reproduced and adapted from \cite{NectarIvan2023} under the Creative Commons Attribution 3.0 license.}
	\label{fig:elastic_film_shematic}
\end{figure}

To this end, consider a Newtonian fluid in a closed trough with horizontal length $2\ell$, as shown in Fig.~\ref{fig:elastic_film_shematic}. The fluid has equilibrium height $h_0 \ll \ell$, making it a ``thin'' (or ``slender'') film, such that the lubrication approximation applies \cite{Stone2017LH}. Linearizing about the flat film state, such that $h(x,t)=h_0\cdot[1+u(x,t)]$ with $\max_{x,t}|u(x,t)|\ll1$, and by standard methods \cite{Oron1997}, one obtains \cite{NectarIvan2023} the (long-wave) \emph{thin film equation}:
\begin{equation}
	\frac{\partial u}{\partial t} = \frac{h_0^3}{12\mu_f} \frac{\partial}{\partial x} \left( \rho_f g  \frac{\partial u}{\partial x} + B \frac{\partial^5 u}{\partial x^5} \right) .
	\label{eq:elastic_film_0}
\end{equation}
For the boundary conditions (BCs) considered here (and depicted in Fig.~\ref{fig:elastic_film_shematic}), the film is free to move vertically at the edges of the trough but maintains a $90^\circ$ contact angle (``non-wetting'' condition, $\partial u/\partial x =0$ at $x=\pm\ell$), which ensures that $h=const.$ is an equilibrium state, as in \cite{Liu2023}. It is assumed that the equilibrium state can be accommodated on the domain, precluding buckling \cite{Kodio2017}. At the same time, the sheet should be considered as (at least slightly) compressible, which is typical for sufficiently thin sheets made from polymeric materials \cite{Chandler2020}. 

Further, we enforce no moment at the walls for this non-pinned elastic interface ($\partial^2 u/\partial x^2 =0$ at $x=\pm\ell$). No fluid flux through the lateral confinement requires that $\partial^5 u/\partial x^5 = 0$ at $x=\pm\ell$; recall the flux is related to the quantity in the parentheses on the right-hand side of Eq.~\eqref{eq:elastic_film_0}. The boundary conditions are summarized mathematically in Eq.~\eqref{eq:BC_6} below. Importantly, in \cite{NectarIvan2023}, we showed that the BCs considered herein lead to a \emph{self-adjoint} sixth-order eigenvalue problem associated with Eq.~\eqref{eq:elastic_film_0}. The corresponding eigenfunctions form a complete orthonormal set, a fact that we exploit in this work. Other boundary conditions are also possible, such as the clamped--free ones employed by \citet{Hosoi2004}, which do not enforce no-flux. In that case, imposing clamping and force-free, but not torque-free, conditions ($u=0$, $\partial u/\partial x=0$, and $\partial^3 u/\partial x^3 = 0$ at $x=\pm\ell$) would also lead to a self-adjoint problem, which has not been explored yet. Non-self-adjoint problems also arise and are of interest in their own right \cite{Gabay2023}.

To proceed, we make Eq.~\eqref{eq:elastic_film_0}  dimensionless by introducing the ``hat'' variables via the transformations:
\begin{equation}
	x = \ell\hat{x}, \qquad t = t_c \hat{t},\qquad u(x,t) = u_c \hat{u}(\hat{x},\hat{t}).
	\label{eq:elastic_film_ndvar}
\end{equation}
Substituting the variables from Eq.~\eqref{eq:elastic_film_ndvar} into Eq.~\eqref{eq:elastic_film_0}, we obtain
\begin{equation}
	\frac{\partial \hat{u}}{\partial \hat{t}} = \frac{h_0^3t_c}{12\mu_f} \frac{\rho_f g}{\ell^2}  \frac{\partial^2 \hat{u}}{\partial \hat{x}^2} + \frac{h_0^3t_c}{12\mu_f} \frac{B}{\ell^6} \frac{\partial^6 \hat{u}}{\partial \hat{x}^6}.
	\label{eq:elastic_film}
\end{equation}
Notice that $u_c$ is arbitrary, as it cancels out of the problem due to linearity in $u$.
In \cite{NectarIvan2023}, we focused on the case in which the dominant force on the interface is the elastic bending resistance. Thus, the choice of time scale in \cite{NectarIvan2023} was the viscous--elastic \cite{Elbaz2014} one, $t_{c,ve} = {12\mu_f \ell^6}/{Bh_0^3}$, as in \cite{Tulchinsky2016,MartinezCalvo2020}. This choice makes the coefficient of $\partial^6 u/\partial \hat{x}^6$ in Eq.~\eqref{eq:elastic_film} unity. Likewise, we could also pick the time scale to be the viscous--gravity one, $t_{c,vg}={12\mu_f \ell^2}/{\rho_f g h_0^3}$, which would render the coefficient of $\partial^2 u/\partial \hat{x}^2$ in Eq.~\eqref{eq:elastic_film} unity. The ratio of these two time scales is 
\begin{equation}
	\Bo = \frac{t_{c,ve}}{t_{c,vg}} = \frac{\rho_f g \ell^4}{B},
\end{equation}	
which is termed the \emph{elastic Bond number} (see \cite{Duprat2011,Saeter2024}). Without loss of generality, we may take $t_c = t_{c,ve}$, making $\Bo$ the single controlling parameter of the problem. In particular, we prefer the choice $t_c = t_{c,ve}$ so that we can take $\Bo\to0$ in the governing equations and focus on the ``pure'' sixth-order problem in several examples below. 

We can estimate the value of $\Bo$ for a typical thin film based on prior experiments. \citet{Pedersen2019} placed a $160~\si{\nano\meter}$ thick sheet of polysulfone (with a reported bending rigidity $B \simeq 1.3~\si{\pico\joule}$) on top of a thin fluid film. \citet{Gabay2023} considered a closed trough of axial width $2\ell = 4~\si{\milli\meter}$ filled with a silicone oil of density $\rho_f = 970~\si{\kilo\gram\per\meter\cubed}$. Combining these numbers, we have the estimate $\Bo \simeq 10^5$. Note, however, that $\ell^4$ ``controls'' the magnitude of $\Bo$ once the fluid and elastic interface properties are set; halving the trough length decreases $\Bo$ by a factor of $(0.5)^4=0.065$. Indeed, for the most recent experiments by \citet{Saeter2024}, $\Bo \simeq 59$ based on $B \simeq 5~\si{\milli\joule}$ for a $5~\si{\milli\meter}$ thick sheet of silicon-based elastomer, again on top of a silicon oil film with $\rho_f = 970~\si{\kilo\gram\per\meter\cubed}$ and $\ell \simeq 75~\si{\milli\meter}$. In the latter experiment, the model is posed on an infinite domain, and the scale $\ell$ is taken to represent the support of a compact initial condition. Clearly, a wide range of \emph{large} values of $\Bo$ can be justified depending on the problem at hand, and we discuss the effect of this parameter below.

To this end, in this work, we consider the following dimensionless sixth-order initial boundary value problem (IBVP) for $u=u(x,t)$:
\begin{subequations}
	\begin{empheq}[left = \empheqlbrace]{alignat=2}
		\frac{\partial u}{\partial t} - \Bo \frac{\partial^2 u}{\partial x^2} - \frac{\partial^6 u}{\partial x^6} &= 0, &\qquad -1<x<+1,\quad t>0, \label{eq:PDE_6}\\
		\left.\frac{\partial u}{\partial x}\right|_{x=\pm1} = \left.\frac{\partial^2 u}{\partial x^2}\right|_{x=\pm1} = \left.\frac{\partial^5 u}{\partial x^5}\right|_{x=\pm1} &= 0,\quad &t>0, \label{eq:BC_6}\\
		u(x,0) &= u^0(x), \quad &-1<x<+1, \label{eq:IC_6}
	\end{empheq}\label{eq:IBVP_6}%
\end{subequations}
where $u^0(x)$ is a given function. The remainder of this paper is organized as follows. In Sec.~\ref{sec:green}, we explicitly construct the Green's function for IBVP~\eqref{eq:IBVP_6} with $\Bo=0$ on both a finite domain (using the sixth-order eigenfunctions) and on an infinite domain (using the Fourier transform). Then, in Sec.~\ref{sec:ss_solutions}, we use both the machinery from Sec.~\ref{sec:green} and a Galerkin spectral expansion to derive exact and numerical solutions, respectively, to the finite-interval problem for two representative cases (a point-load and a stepped initial disturbance on the film). Furthermore, in Sec.~\ref{sec:ss_solutions}, we discuss the intermediate asymptotics of the solution for these two model initial conditions, guided by the infinite-line problem's self-similar solutions. We also determine the effect of the elastic Bond number $\Bo$ on the typical adjustment time of the film.  Conclusions and avenues for future work are stated in Sec.~\ref{sec:conclusion}. The expressions for the sixth-order eigenfunctions and key details from \cite{NectarIvan2023} are summarized in Appendix~\ref{sec:CON} for completeness.

\section{Green's functions for $\Bo=0$}\label{sec:green}

To solve IBVP~\eqref{eq:IBVP_6} with $\Bo=0$, we will employ the complete orthonormal set of eigenfunctions derived in \cite{NectarIvan2023} and summarized in Appendix~\ref{sec:CON}. In Sec.~\ref{sec:green_fin}, we will construct the Green's function explicitly. Then, in Sec.~\ref{sec:green_inf}, we will review the construction of the Green's function of IBVP~\eqref{eq:IBVP_6} with $\Bo=0$ on the real line, which was also discussed in the Appendix of \cite{Pedersen2021}.

\subsection{1D Green's function for the problem on a finite interval}\label{sec:green_fin}

Following \citet[Sec.~5.2]{Duffy2015}, the finite-interval Green's function, $G_\mathrm{fin}(x,t|\xi,\tau)$, of the operator in Eq.~\eqref{eq:PDE_6} with $\Bo=0$ satisfies
\begin{subequations}
	\begin{empheq} [left = \empheqlbrace]{alignat=2}
		\frac{\partial G_\mathrm{fin}}{\partial t} - \frac{\partial^6 G_\mathrm{fin}}{\partial x^6} &= \delta(x-\xi)\delta(t-\tau), \qquad &-1<x,\xi<+1,\quad t,\tau>0, \label{eq:Green_PDE_6}\\
		\left.\frac{\partial G_\mathrm{fin}}{\partial x}\right|_{x=\pm1} = \left.\frac{\partial^2 G_\mathrm{fin}}{\partial x^2}\right|_{x=\pm1} &= \left.\frac{\partial^5 G_\mathrm{fin}}{\partial x^5}\right|_{x=\pm1} = 0,\quad &t>0,\\
		G_\mathrm{fin}(x,0|\xi,\tau) &= 0, \quad &-1<x<+1,\label{eq:Green_IC_6}
	\end{empheq}\label{eq:Green_IBVP_6}%
\end{subequations}
where $\delta(\,\cdot\,)$ is Dirac's delta distribution. Taking a Laplace transform, such that $\bar{G}_\mathrm{fin}(x,s|\xi,\tau) = \int_0^\infty \re^{-st} G_\mathrm{fin}(x,t|\xi,\tau) \,\rd t$, of Eq.~\eqref{eq:Green_PDE_6}, and using the initial condition~\eqref{eq:Green_IC_6}, we obtain
\begin{subequations}\label{eq:Green_BVP_6}
\begin{equation}
	s \bar{G}_\mathrm{fin} - \frac{\rd^6 \bar{G}_\mathrm{fin}}{\rd x^6} = \delta(x-\xi)\re^{-s\tau},\qquad -1<x<+1,
\end{equation}
subject to
\begin{equation}
	\left.\frac{\partial \bar{G}_\mathrm{fin}}{\partial x}\right|_{x=\pm1} = \left.\frac{\partial^2 \bar{G}_\mathrm{fin}}{\partial x^2}\right|_{x=\pm1} = \left.\frac{\partial^5 \bar{G}_\mathrm{fin}}{\partial x^5}\right|_{x=\pm1} = 0,\qquad t>0. \label{eq:Green_BC_6}
\end{equation}
\end{subequations}

BVP~\eqref{eq:Green_BVP_6} gives rise to an associated eigenvalue problem (EVP):
\begin{subequations}
	\label{eq:main_6th_EVP}
	\begin{align}
		-\frac{\rd^6 \psi_n}{\rd x^6} &= \lambda_n^6 \psi_n , 	\\
		\left.\frac{\rd \psi_n}{\rd x}\right|_{x=\pm1} = \left.\frac{\rd^2 \psi_n}{\rd x^2}\right|_{x=\pm1}&= \left.\frac{\rd^5 \psi_n}{\rd x^5}\right|_{x=\pm1}=0 .
	\end{align}
\end{subequations}	
The eigenfunctions solving EVP~\eqref{eq:main_6th_EVP} were explicitly constructed and shown to form a complete orthonormal set on $L^2[-1,+1]$ in \cite{NectarIvan2023}. Thus, the solution of BVP~\eqref{eq:Green_BVP_6} can be written as
\begin{equation}
	\bar{G}_\mathrm{fin}(x,s|\xi,\tau) = \re^{-s\tau} \sum_{n=0}^\infty \frac{\psi_n(\xi)}{s +\lambda_n^6}\psi_n(x).
	\label{eq:Green_2D_interval_Laplace_domain}
\end{equation}
Inverting the Laplace transform using tables of inverses \cite{BoyceDiPrima_10thEd}, Eq.~\eqref{eq:Green_2D_interval_Laplace_domain} becomes:
\begin{equation}
	G_\mathrm{fin}(x,t|\xi,\tau) = H(t-\tau) \left[ \sum_{n=0}^\infty \re^{-\lambda_n^6(t-\tau)} \psi_n(\xi) \psi_n(x)  \right],
	\label{eq:Green_2D_interval}
\end{equation}
where $H(\,\cdot\,)$ is Heaviside's unit step function. 

Using the Green's function, one can write down the general solution to IBVP~\eqref{eq:IBVP_6} with $\Bo=0$ as
\begin{equation}
	u(x,t) = \int_{-1}^{+1} G_\mathrm{fin}(x,t|\xi,0) u^0(\xi) \, \rd\xi.
\end{equation}
Introducing Eq.~\eqref{eq:Green_2D_interval} into the latter and defining $u^0_n = \int_{-1}^{+1} \psi_n(\xi) u^0(\xi) \, \rd\xi$ as the expansion coefficients of the initial condition, we have the general result for $\Bo=0$:
\begin{equation}
	u(x,t) = H(t)\left[ \sum_{n=0}^\infty \re^{-\lambda_n^6t} u^0_n \psi_n(x)  \right].
	\label{eq:Green_2D_interval_convl}
\end{equation}
Note that Eq.~\eqref{eq:Green_2D_interval_convl} is the same expression as the one that would be obtained by solving IBVP~\eqref{eq:IBVP_6} (with $\Bo=0$) by the separation of variables method.

So far, we have denoted the set of eigenfunctions simply as $\{ \psi_n(x) \}_{n=0,1,2\hdots}$. Below, we will split the eigenfunctions into odd (``sine'') $\{ \psi_n^s(x) \}_{n=1,2\hdots}$ and even (``cosine'')  $\{ \psi_n^c(x) \}_{n=0,1,2\hdots}$, as discussed in \cite{NectarIvan2023}. This choice will make the solution of model problems simpler but will require special treatment of the $n=0$ term (see Eq.~\eqref{eq:delta_exp_defn} and Eq.~\eqref{eq:u_x_expansion_formula}).

\subsection{1D Green's function for the problem on the real line}\label{sec:green_inf}

Now, consider an infinite domain, $-\infty<x<+\infty$. This case was also treated in the Appendix of \citet{Pedersen2021} (see also the related examples in \cite{Tulchinsky2016,Rubin2017}), but we include it here for completeness and comparison to the previous result. Again, following \citet[Sec.~5.1]{Duffy2015}, the Green's function, $G_\mathrm{inf}(x,t|\xi,\tau)$, of the operator in Eq.~\eqref{eq:PDE_6} with $\Bo=0$ now satisfies
\begin{subequations}
	\begin{empheq} [left = \empheqlbrace]{alignat=2}
		\frac{\partial G_\mathrm{inf}}{\partial t} - \frac{\partial^6 G_\mathrm{inf}}{\partial x^6} &= \delta(x-\xi)\delta(t-\tau), \qquad &-\infty<x,\xi<+\infty,\quad t,\tau>0, \label{eq:Green_PDE_6_inf}\\
		\lim_{|x|\to\infty} G_\mathrm{inf}(x,t|\xi,\tau) &= 0,\quad &t>0,\\
		G_\mathrm{inf}(x,0|\xi,\tau) &= 0, \quad &-\infty<x<+\infty.\label{eq:Green_IC_6_inf}
	\end{empheq}
\end{subequations}
Again, we first take the Laplace transform of Eq.~\eqref{eq:Green_PDE_6_inf} and use the initial condition~\eqref{eq:Green_IC_6_inf}, to obtain
\begin{equation}
	s \bar{G}_\mathrm{inf} - \frac{\rd^6 \bar{G}_\mathrm{inf}}{\rd x^6} = \delta(x-\xi)\re^{-s\tau},\qquad -\infty<x<+\infty.
\end{equation}
Now, we take the Fourier transform, such that $\hat{\bar{G}}_\mathrm{inf}(k,s|\xi,\tau) = \int_{-\infty}^{+\infty} \re^{-\ri k x} \bar{G}_\mathrm{inf}(x,s|\xi,\tau) \,\rd x$, of the last equation to obtain
\begin{equation}
	\hat{\bar{G}}_\mathrm{inf}(k,s|\xi,\tau)  = \frac{\re^{-\ri k \xi}\re^{-s\tau}}{s + k^6}.
	\label{eq:Green_2D_line_dual_transform}
\end{equation}

Inverting the Laplace transform using tables of inverses \cite{BoyceDiPrima_10thEd}, Eq.~\eqref{eq:Green_2D_line_dual_transform} becomes:
\begin{equation}
	\hat{G}_\mathrm{inf}(k,t|\xi,\tau) = H(t-\tau) \re^{-k^6(t-\tau) - \ri k\xi}.
\end{equation}
The Fourier inversion integral, $G_\mathrm{inf}(x,t|\xi,\tau) =  \frac{1}{2\pi}\int_{-\infty}^{+\infty} \re^{\ri k x} \hat{G}_\mathrm{inf}(k,t|\xi,\tau) \,\rd k$, is challenging. Using \textsc{Mathematica}~\cite{Mathematica}, we find
\begin{subequations}
	\label{eq:Green_2D_line_all}
	\begin{align}
	\label{eq:Green_2D_line}
	G_\mathrm{inf}(x,t|\xi,\tau) &= \frac{H(t-\tau)}{(t-\tau)^{1/6}} \mathcal{G}(\zeta), 	\qquad \zeta:= \frac{x-\xi}{(t-\tau)^{1/6}},\\[2mm]
	\label{eq:attractor_delta}
	\mathcal{G}(\zeta) &:= 
		\frac{\Gamma(7/6)}{\pi}  ~{}_0F_4\left[\left\{\frac{1}{3},\frac{1}{2},\frac{2}{3},\frac{5}{6}\right\}, -\left(\frac{\zeta}{6}\right)^6 \right] 
		- \zeta^2\frac{1}{12\sqrt{\pi}} ~{}_0F_4\left[\left\{\frac{2}{3},\frac{5}{6},\frac{7}{6},\frac{4}{3}\right\},  - \left(\frac{\zeta}{6}\right)^6 \right] \nonumber\\
		&\phantom{:=} + \zeta^4\frac{1}{432\,\Gamma(7/6)} ~{}_0F_4\left[\left\{\frac{7}{6},\frac{4}{3},\frac{3}{2},\frac{5}{3}\right\}, -\left(\frac{\zeta}{6}\right)^6 \right],
	\end{align}	
\end{subequations}
where $\Gamma(\cdot)$ is Euler's Gamma function, and ${}_0F_4(\{\cdot,\cdot,\cdot,\cdot\},z)$ denotes the generalized hypergeometric function ${}_pF_q\left(\substack{\mathbf{a}=\{a_1,\hdots,a_p\}\\ \mathbf{b}=\{b_1,\hdots,b_q\}};z\right)$ with $p=0$ ($\mathbf{a}=\mathbf{0}$) and $q=4$ ($\mathbf{b}=\{\cdot,\cdot,\cdot,\cdot\}$) \cite[Ch.~16]{NIST:DLMF}. 
Finally, the general solution to the version of IBVP~\eqref{eq:IBVP_6} (with $\Bo=0$) on the real line is
\begin{equation}
	u(x,t) = \int_{-\infty}^{+\infty} G_\mathrm{inf}(x,t|\xi,0) u^0(\xi) \, \rd\xi.
	\label{eq:Green_2D_line_convl}
\end{equation}
Note that while $\zeta$ in Eq.~\eqref{eq:Green_2D_line} was defined in terms of $x$, $\xi$, $t$ and $\tau$ in the context of the Green's function, below we shall use the same notation even if $\xi,\tau=0$, in the context of the solution, say $u(x,t)=u(\zeta)$ with $\zeta := x/t^{1/6}$, without fear of confusion. The context will unambiguously indicate the meaning of $\zeta$.

In closing this subsection, we note that $\mathcal{G}(\zeta)$ also solves $6\mathcal{G}''''' + \zeta \mathcal{G} = 0$, where $'=\rd/\rd \zeta$, with suitable decay conditions at infinity. This ordinary differential equation (ODE) can be obtained by assuming a self-similar solution $u(x,t)=\mathcal{G}(\zeta)/t^{1/6}$, $\zeta = x/t^{1/6}$ and substituting this ansatz into Eq.~\eqref{eq:IBVP_6} with $\Bo=0$. However, such a high-order ODE with a nonconstant coefficient must again be solved by \textsc{Mathematica} or by referring to tables and books of special functions, yielding the result in Eq.~\eqref{eq:Green_2D_line_all} after solving for the integration constants.

\section{Self-similar asymptotics of the finite-interval solution}\label{sec:ss_solutions}

In this section, we will solve the IBVP~\eqref{eq:IBVP_6} for two representative choices of $u^0(x)$. First, in Sec.~\ref{sec:ss_exact} for $\Bo=0$, we will employ the Green's function derived in Sec.~\ref{sec:green_fin}. Next, in Sect.~\ref{sec:ss_numerical}, we will introduce a Galerkin spectral expansion based on the same sixth-order eigenfunctions used in deriving the Green's function.

\subsection{Problems with exact solutions}\label{sec:ss_exact}

\subsubsection{Point-load initial condition}\label{sec:ss_delta}

First, we summarize the expansion of the Dirac $\delta$ into the basis of sixth-order eigenfunctions, which will furnish the coefficients $u_n^0$ in Eq.~\eqref{eq:Green_2D_interval_convl}. The Dirac $\delta$ can be considered as an even ``function''. Hence, only the even eigenfunctions $\psi^c_n(x)$ (including $\psi^c_0(x)=1$) appear in its spectral expansion into the chosen basis of sixth-order eigenfunctions:
\begin{subequations}
	\label{eq:delta_exp} 
	\begin{alignat}{3}
 		\delta(x) &= \frac{1}{2}\delta^c_0\,\psi_0^c(x) + \sum_{n=1}^{\infty} \delta^c_n\,\psi_n^c(x), & \label{eq:delta_exp_defn}\\
  	    \delta^c_0 &= 1,  &&\qquad \qquad \qquad (n=0) \\
		\delta^c_n &= c_n^c \times \Bigg[ \frac{4 \sin{\frac{\lambda_n^c}{2}} \sin{\lambda_n^c} \cosh{\frac{\sqrt{3} \lambda_n^c}{2}}}{\cos{\lambda_n^c}-\cosh{\sqrt{3}\lambda_n^c}}+1\Bigg],  &&\qquad \qquad \qquad (n\ge1)
	\end{alignat}
\end{subequations}
where the expression for $c_n^c$ is given by Eqs.~\eqref{eq:c_n^c_defn} and \eqref{eq:d_n^c_defn}.

Observe that, in this case, the solution~\eqref{eq:Green_2D_interval_convl} takes the particularly simple form $u(x,t) = G_\mathrm{fin}(x,t|0,0)$, being given directly by the Green's function. Similarly, the infinite line solution is simply $u(x,t) = G_\mathrm{inf}(x,t|0,0) = \mathcal{G}(\zeta)$, which is a self-similar function of $x$ and $t$, as seen from Eq.~\eqref{eq:attractor_delta}. Importantly, $\mathcal{G}(\zeta)$ is termed the self-similar intermediate \emph{attractor} \cite{Barenblatt1972,Barenblatt1996}, which we expect also to represent the rescaled finite-interval solution $G_\mathrm{fin}(x,t|0,0)$ (which is \emph{not} explicitly a self-similar function of $x$ and $t$), at least until the effects of the finite boundaries are felt, for any initial conditions with finite zeroth moment.

Note that taking $u^0(x) = \delta(x)$ is, in a sense, the perfect or ideal case in which we would expect the solution to lie on the attractor right away from $t=0$. Starting with nonideal initial conditions, such as a tall and narrow post that approximates $\delta(x)$, would more clearly highlight the transition \emph{onto} the attractor for intermediate times, and \emph{away} from it for late times. Such an example is shown in Appendix~\ref{sec:narrow_support}. We focus our main discussion below on the ideal case(s).

\subsubsection{Stepped initial condition}\label{sec:ss_heaviside}

Now, consider a ``stepped'' initial condition, motivated by previous theoretical and experimental work on the leveling of capillary thin films \cite{Salez2012,McGraw2012}. We represent this type of initial condition by the shifted step function
\begin{equation}
	\label{eq:odd_step_defn}
	H_\mathrm{odd}(x) :=  \frac{1}{2}\sgn(x)
	= 
	\begin{cases}
	 	+\frac{1}{2}, &\quad x \leq 0 ,\\[1mm]
		-\frac{1}{2}, &\quad x > 0 ,
	\end{cases}
\end{equation}
where $\sgn(x) = x/|x|$ for $x\ne0$ while $\sgn(0)=0$ is the signum function. 
We choose to define an odd step function as the initial condition so that only the odd eigenfunctions $\psi^s_n(x)$ are required to obtain its spectral expansion, which leads to computational efficiency (by requiring fewer terms in the series of $H_\mathrm{odd}$ compared to the Heaviside step function $H$) as well as conceptual simplicity. 

For this case, the coefficients $u_n^0$ in Eq.~\eqref{eq:Green_2D_interval_convl} will come from the expansion of $H_\mathrm{odd}(x)$ into the basis of sixth-order eigenfunctions, which reads as follows:
\begin{subequations}
	\label{eq:odd_step_exp} 
	\begin{align}
		H_\mathrm{odd}(x)  &= \sum_{n=1}^{\infty} h^s_n\,\psi_n^s(x), \\
		h^s_n &= \frac{c_n^s}{\lambda^s_n  \left(\cos{\lambda^s_n} + \cosh{\sqrt{3} \lambda^s_n}\right)} \times \Bigg[-\cos{\lambda^s_n} + \left(\cos{\frac{\lambda^s_n}{2}} + \cos{ \frac{3 \lambda^s_n}{2}}\right) \cosh{\frac{\sqrt{3} \lambda^s_n}{2}} \\
		&\phantom{=}-\cosh{\sqrt{3} \lambda^s_n} - \sqrt{3}  \left(\sin{\frac{\lambda^s_n}{2}} - \sin{\frac{3 \lambda^s_n}{2}}\right) \sinh{\frac{\sqrt{3} \lambda^s_n}{2}} \Bigg], \nonumber
	\end{align}%
\end{subequations}
where the expression for $c_n^s$ is given in Eqs.~\eqref{eq:c_n^s_defn} and \eqref{eq:d_n^s_defn}.

The corresponding infinite-line solution and expected self-similar attractor can be found from Eq.~\eqref{eq:Green_2D_line_convl} upon substituting Eq.~\eqref{eq:odd_step_defn} for $u^0(\xi)$. Performing the laborious integration with the help of \textsc{Mathematica} and using the properties of ${}_pF_q$, we obtain the self-similar expression:
\begin{multline}
	u(x,t) =  \mathcal{H}(\zeta) := - \frac{\zeta}{6\pi} \Bigg\{ 
	\Gamma(1/6) ~{}_1F_5\left[\substack{\{\frac{1}{6}\}\\[1mm] \{\frac{1}{3},\frac{1}{2},\frac{2}{3},\frac{5}{6},\frac{7}{6}\}}; -\left(\frac{\zeta}{6}\right)^6 \right] 
	- \zeta^2 \frac{\sqrt{\pi}}{6} ~{}_1F_5\left[\substack{\{\frac{1}{2}\}\\[1mm] \{\frac{2}{3},\frac{5}{6},\frac{7}{6},\frac{4}{3},\frac{3}{2}\}}; -\left(\frac{\zeta}{6}\right)^6 \right] \\
	+ \zeta^4 \frac{\Gamma(5/6)}{120} ~{}_1F_5\left[\substack{\{\frac{5}{6}\}\\[1mm] \{\frac{7}{6},\frac{4}{3},\frac{3}{2},\frac{5}{3},\frac{11}{6}\}}; -\left(\frac{\zeta}{6}\right)^6 \right] \Bigg\}, 
	\qquad \zeta := \frac{x}{t^{1/6}}.
	\label{eq:attractor_2D_line_odd_step}
\end{multline}
Again, we expect that $\mathcal{H}(\zeta)$ represents the rescaled finite-interval solution (which is \emph{not} explicitly a self-similar function of $x$ and $t$), at least until the effect of the finite boundaries is felt.

Note that $H_\mathrm{odd}(x)$ is a step of ``height'' 1 at $x=0$ and resembles the initial condition used in one of the cases presented in \cite{Pedersen2021}, allowing for a further qualitative comparison. Thus, as expected, Eq.~\eqref{eq:attractor_2D_line_odd_step} is also similar to Eq.~(A\,7) in \cite{Pedersen2021}, which is the expression corresponding to taking $u^0(\xi) = H(-x)$ (see also the related examples in \cite{Tulchinsky2016,Rubin2017}). As in the previous subsection, this is an ideal case in which we expect the solution to lie on the attractor right away from $t=0$.

\subsection{Numerical approach to IBVP~\eqref{eq:IBVP_6}}\label{sec:ss_numerical}

For the case $\Bo \ne 0$, we solve IBVP~\eqref{eq:IBVP_6} numerically. The reason is that both the Green's function approach followed in Sec.~\ref{sec:green_fin} and the separation of variables method would both require the solution of a different eigenvalue problem (other than the one given in Eq.~\eqref{eq:main_6th_EVP}), containing a second derivative term, and the derivation of a new complete set of orthonormal eigenfunctions. Here, we would like to demonstrate that this rather toilsome process is unnecessary once an eigenfunction expansion for the highest-derivative operator is obtained. Specifically, we showed in \cite{NectarIvan2023} that a classical Galerkin spectral method based on the (even and odd sets of) eigenfunctions that solve Eq.~\eqref{eq:main_6th_EVP} can efficiently and accurately solve sixth-order equations containing lower (second and/or fourth) order derivative terms, providing the BVP is still subject to the same BCs~\eqref{eq:BC_6}.

\subsubsection{Semi-discretization scheme}

To begin, we expand the sought function $u(x,t)$ solving IBVP~\eqref{eq:IBVP_6} in a series in terms of the chosen sixth-order eigenfunctions as:
\begin{equation}
	u(x,t) \approx  \frac{1}{2}\hat{u}_0^c \psi_0^c(x) + \sum_{n=1}^{M} \hat{u}_n^c(t) \psi_n^c(x) + \sum_{n=1}^{M} \hat{u}_n^s(t) \psi_n^s(x).
	\label{eq:u_x_expansion_formula}
\end{equation}
The expansion~\eqref{eq:u_x_expansion_formula} is introduced into Eq.~\eqref{eq:PDE_6}. Then, taking successive inner products with the $\psi_{m}^c(x)$ and $\psi_{m}^s(x)$, for $m=1,2,\ldots,M$, we obtain the following  system of linear ordinary differential equations (ODEs) for the time-dependent spectral coefficients $\hat{u}_{m}^c(t)$ and $\hat{u}_{m}^s(t)$:
\begin{subequations}
	\label{eq:gen_DS_DEs} 
	\begin{align}
		\frac{\dif \hat{u}_{m}^c(t)}{\dif t}& =  \sum_{n=1}^M  \left[\Bo \beta_{nm}^c - (\lambda^c_{m})^6 \delta_{nm} \right]\hat{u}_n^c(t) , 
		\qquad (m=1,2,\ldots,M) 	\label{eq:gen_DS_even_DEs} \\
		\frac{\dif \hat{u}_{m}^s(t)}{\dif t}&= \sum_{n=1}^M\left[\Bo  \beta_{nm}^s  - (\lambda^s_{m})^6 \delta_{nm}\right]\hat{u}_{n}^s(t) , 
		\qquad (m=1,2,\ldots,M)	\label{eq:gen_DS_odd_DEs} 
	\end{align}
\end{subequations}	
where $\delta_{nm}$ is Kronecker's delta, and we have defined 
\begin{subequations}
		\label{eq:sec_der_def}
\begin{align}
	\beta_{nm}^c &:= \int_{-1}^{+1} (\psi_n^c)''  \psi_m^c \,\rd x, \qquad (n,m=1,2,\ldots,M) 	\label{eq:sec_der_even_def}\\
	\beta_{nm}^s &:= \int_{-1}^{+1} (\psi_n^s)''  \psi_m^s \,\rd x, \qquad (n,m=1,2,\ldots,M) 	\label{eq:sec_der_odd_def}
\end{align}
\end{subequations}
as the expansions of the second derivatives of the eigenfunctions back into the basis, with $x$-derivatives denoted by primes. The expressions for the coefficients $\beta_{nm}^c$ and $\beta_{nm}^s$ are given in Appendix~\ref{sec:CON} by Eq.~\eqref{eq:sec_deriv_even_form} and  Eq.~\eqref{eq:sec_deriv_odd_form} respectively.

The coefficient $\hat{u}_0^c$ in the expansion~\eqref{eq:u_x_expansion_formula} is not time-dependent and is found separately from the initial condition~\eqref{eq:IC_6}:
\begin{equation}
	\hat{u}_0^c = \int_{-1}^{+1} u^0(x)\dif x .
\end{equation}	
The ICs for the dynamical system~\eqref{eq:gen_DS_DEs} are the remaining coefficients in the spectral expansion of $u^0(x)$, namely:
 \begin{subequations}
 	\label{eq:gen_DS_ICs} 
 	\begin{align}
         \hat{u}_m^{c}(0) := \hat{u}_m^{0,c}&= \int_{-1}^{+1}  u^0(x)   \psi_m^c(x) \,\rd x, \qquad (m=1,2,\ldots,M)\label{eq:gen_DS_ICs_even} \\
         \hat{u}_m^{s}(0) :=  \hat{u}_m^{0,s}&= \int_{-1}^{+1}  u^0(x)   \psi_m^s(x) \,\rd x. \qquad (m=1,2,\ldots,M)\label{eq:gen_DS_ICs_odd}
 	\end{align}	
\end{subequations}	
We mention here that the dynamical system~\eqref{eq:gen_DS_even_DEs} subject to the  ICs~\eqref{eq:gen_DS_ICs_even} for the even coefficients $\{\hat{u}_m^c(t)\}_{m=1}^M$ is decoupled from the dynamical system~\eqref{eq:gen_DS_odd_DEs} subject to the ICs~\eqref{eq:gen_DS_ICs_odd} for the odd coefficients $\{\hat{u}_{m}^s(t)\}_{m=1}^M$. Therefore, we can solve for the odd and even coefficients as separate linear systems.

Employing a matrix-vector formulation, we denote the sought vectors of coefficients by 
\begin{subequations}
	\begin{align}
		\bm{\hat{u}^c}(t) &= \big[ \hat{u}_{1}^c(t),\hat{u}_{2}^c(t),\hdots,\hat{u}_{M}^c(t) \big]^{\mathrm{T}}, \\
		\bm{\hat{u}^s}(t) &= \big[ \hat{u}_{1}^s(t),\hat{u}_{2}^s(t),\hdots,\hat{u}_{M}^s(t) \big]^{\mathrm{T}},
	\end{align}%
\end{subequations}	
and rewrite the dynamical systems~\eqref{eq:gen_DS_DEs} as
\begin{subequations} 
	 \label{eq:gen_Matr_DS_DEs}
	\begin{align}
		\frac{\dif \bm{\hat{u}^c}}{\dif t}&= \bm{\mathrm{A}^c} \bm{\hat{u}^c} , 
		\label{eq:gen_Matr_DS_DEs_even} \\
		\frac{\dif \bm{\hat{u}^s}}{\dif t}&= \bm{\mathrm{A}^s} \bm{\hat{u}^s} , 
		\label{eq:gen_Matr_DS_DEs_odd}
	\end{align}%
\end{subequations}	
with the matrices $\bm{\mathrm{A}^c}$ and $\bm{\mathrm{A}^s}$ for the even and odd spectral coefficients, respectively, being given by 
\begin{subequations} 
	\begin{align}
		\bm{\mathrm{A}^c} &= [a_{nm}^c] =  \Bo \beta_{nm}^c  - (\lambda^c_{m})^6 \delta_{nm}, \qquad (n,m = 1,\ldots,M) \label{eq:gen_DS_Matr_defn_even} \\
		\bm{\mathrm{A}^s} &= [a_{nm}^s] =  \Bo \beta_{nm}^s  - (\lambda^s_{m})^6 \delta_{nm}. \qquad (n,m = 1,\ldots,M) \label{eq:gen_DS_Matr_defn_odd} 
	\end{align}%
\end{subequations}
The initial conditions~\eqref{eq:gen_DS_ICs} take the form
\begin{subequations}
	\label{eq:gen_DS_ICs_vec} 
	\begin{align}
		\bm{\hat{u}^{0,c}}&=  \big[ \hat{u}_{1}^{0,c},\hdots,\hat{u}_{M}^{0,c} \big]^{\mathrm{T}}, \label{eq:gen_DS_ICs_even_vec} \\
		\bm{\hat{u}^{0,s}}&=  \big[ \hat{u}_{1}^{0,s},\hdots,\hat{u}_{M}^{0,s} \big]^{\mathrm{T}}. \label{eq:gen_DS_ICs_odd_vec}
	\end{align}%
\end{subequations}

Note that we previously showed \cite{NectarIvan2023} that, for smooth functions, the spectral coefficients $\{\hat{u}_m^c\}_{m=1}^M$  and $\{\hat{u}_m^s\}_{m=1}^M$ decay rapidly with the index $m$. This property of the series allows for highly accurate results to be obtained, even when the expansion~\eqref{eq:u_x_expansion_formula} is terminated at a relatively modest value of $M$. Therefore, even though the matrices $\bm{\mathrm{A}^c}$ and $\bm{\mathrm{A}^s}$ are full, they are relatively small in size.

To integrate the dynamical systems~\eqref{eq:gen_Matr_DS_DEs_even}, \eqref{eq:gen_DS_ICs_even_vec} and \eqref{eq:gen_Matr_DS_DEs_odd}, \eqref{eq:gen_DS_ICs_odd}, we use \textsc{Mathematica}'s highly sophisticated \texttt{NDSolve} function, which offers several advantages, including the ability to control the working precision of the solver (allowing for higher accuracy in its internal computations) and the ability to output the solution as an interpolating function object, enabling the easy evaluation of the solution at any $t\in(0,T]$ \cite{MathematicaDocumentation}. To integrate an IVP, such \eqref{eq:gen_Matr_DS_DEs_even}, \eqref{eq:gen_DS_ICs_even_vec} or \eqref{eq:gen_Matr_DS_DEs_odd}, \eqref{eq:gen_DS_ICs_odd}, \texttt{NDSolve} selects \textsc{Mathematica}'s implementation of LSODA, a versatile adaptive solver that uses a non-stiffness detector to automatically switch between a stiff Gear backward differentiation formula (BDF) method and an Adams non-stiff method \cite{MathematicaDocumentation,Petzold1983}. The suitability of LSODA is discussed in Sec.~\ref{sec:eigenvalues}, as part of a parametric study of the eigenvalues of the matrices $\bm{\mathrm{A}^c}$ and $\bm{\mathrm{A}^s}$.

\subsection{Results and discussion}\label{sec:results}

Using the analytical and numerical methods introduced above, we would now like to discuss the solution behavior for $\Bo \ge 0$ for the two model initial conditions. First, in Sec.~\ref{sec:sss_results_delta}, we discuss the properties of the solution evolving from a point-load initial condition, as represented by the series~\eqref{eq:Green_2D_interval_convl} with coefficients~\eqref{eq:delta_exp}, found using the finite-interval Green's function. Next, in Sec.~\ref{sec:sss_results_step}, we discuss the properties of the solution evolving from a stepped initial condition, as represented by the series~\eqref{eq:Green_2D_interval_convl} with coefficients~\eqref{eq:odd_step_exp}. Finally, in Sec.~\ref{sec:eigenvalues}, we discuss the leading eigenvalue of the spatial operator in Eq.~\eqref{eq:PDE_6}, as a function of $\Bo$, for the example initial conditions, to demonstrate the stability of the Galerkin method and estimate the adjustment time of the film.

\subsubsection{Point-load initial condition}\label{sec:sss_results_delta}

\begin{figure}[ht]
	\centering
	\includegraphics[width=\textwidth]{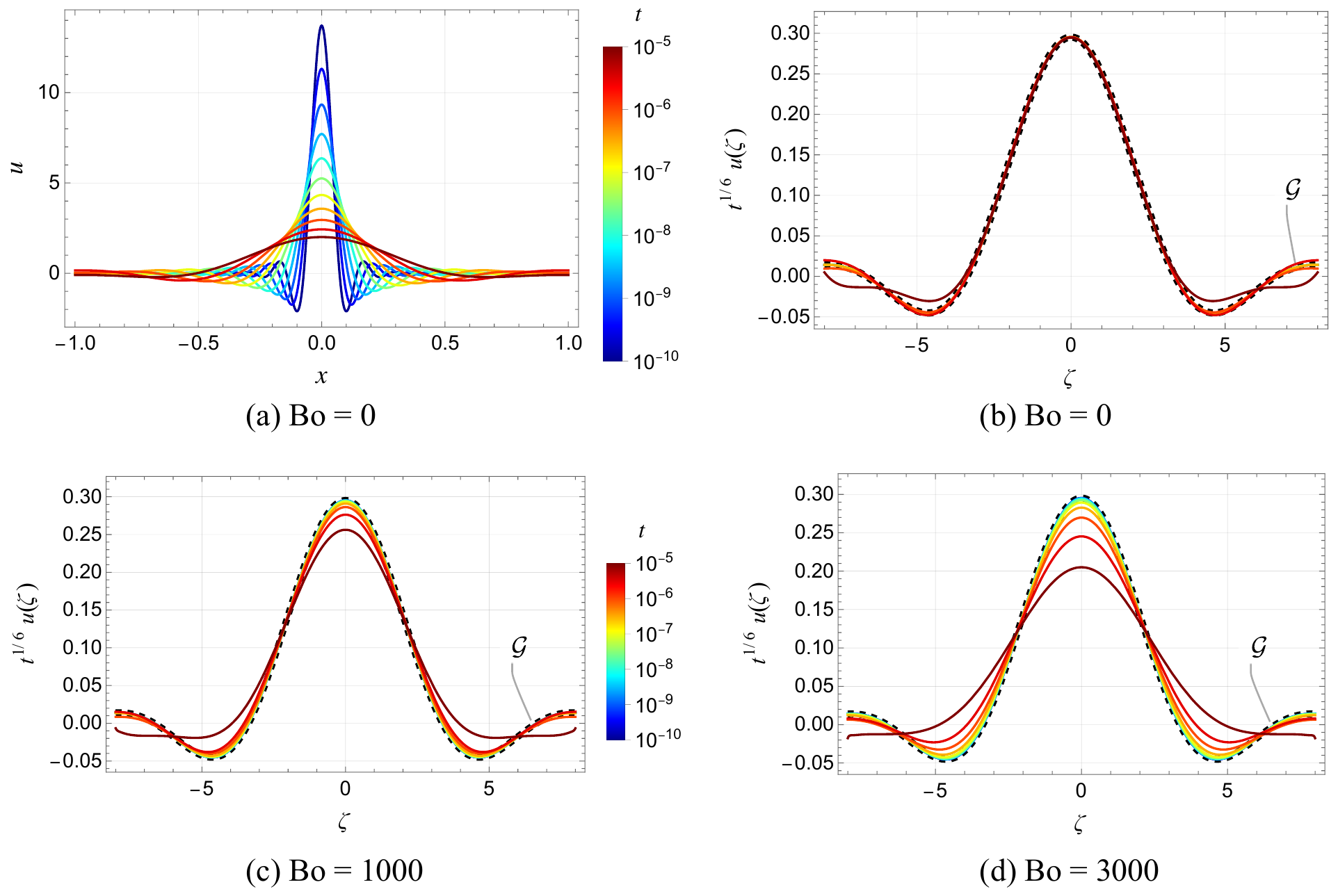}
	\caption{(a) Time-evolution of the solution $u(x,t)$ of IBVP~\eqref{eq:IBVP_6}, given by Eqs.~\eqref{eq:Green_2D_interval_convl} and \eqref{eq:delta_exp}, with $u^0(x) = \delta(x)$ for $\Bo=0$. (b) The profiles from (a) under the self-similar rescaling, i.e., as $t^{1/6} u(\zeta)$ vs.\ $\zeta = x/t^{1/6}$. The thick dashed black curve indicated with $\mathcal{G}$ is the self-similar attractor \eqref{eq:attractor_delta} found via the infinite-line Green's function. (c,d) The equivalent of (b) but for the numerical solution $u(x,t)$ of IBVP~\eqref{eq:IBVP_6} with $u^0(x) = \delta(x)$ for $\Bo=1000$ and $3000$, respectively. For all plots, the spectral series was truncated at $M=50$ (a total of 51 terms).}
\label{fig:delta_function_IC}
\end{figure}

Here, we examine the solution of IBVP~\eqref{eq:IBVP_6} with $u^0(x)=\delta(x)$, for $\Bo=0$. Figure~\ref{fig:delta_function_IC}(a,b) shows the profile of the finite-interval series solution~\eqref{eq:Green_2D_interval_convl} with coefficients~\eqref{eq:delta_exp} for eleven different times $t\in[10^{-10},10^{-5}]$. In Fig.~\ref{fig:delta_function_IC}(a), the solution is plotted in terms of the original dimensionless (unscaled) variables $(x,t)$. We observe that, as expected, the solution decays with time. Also, for the selected time values, the plotted profiles appear localized and are not influenced by the boundaries at early times. Meanwhile, in Fig.~\ref{fig:delta_function_IC}(b), the profiles from (a) are plotted after the self-similar rescaling, i.e., as $t^{1/6} u(\zeta)$ vs.\ $\zeta = x/t^{1/6}$, inferred from the form of the infinite-line Green's function $G_\mathrm{inf}(x,t|0,0)$ in Eqs.~\eqref{eq:Green_2D_line_all}. Thus, in Fig.~\ref{fig:delta_function_IC}(b), the expected self-similar attractor $\mathcal{G}(\zeta)$ from Eq.~\eqref{eq:attractor_delta} is also shown for comparison. 

Figure~\ref{fig:delta_function_IC}(b) portrays the collapse of the rescaled solutions onto a common profile, which demonstrates that the finite-interval series solution~\eqref{eq:Green_2D_interval_convl}, \eqref{eq:delta_exp} does indeed collapse on the self-similar profile, \emph{despite} the series solution not explicitly being a function of a self-similar variable $\zeta$. Importantly, the self-similar profile appears indistinguishable from its infinite-interval counterpart found from Eqs.~\eqref{eq:Green_2D_line_all}. Moreover, as suggested by the previous discussion in Sec.~\ref{sec:ss_delta}, the self-similar solution represents only the \emph{intermediate} asymptotics, and the solution profiles do indeed separate from the attractor as time increases. Eventually, the effect of the finite boundaries comes into play, and the rescaled solution is significantly different from $\mathcal{G}$ for the largest $|\zeta|$ shown. (As noted at the end of Sec.~\ref{sec:ss_delta}, since we start from the ideal $u^0(x)=\delta(x)$ IC, the solution \emph{starts} on the attractor at $t=0$. A parallel nonideal case is presented in Appendix~\ref{sec:narrow_support}.) Beyond the intermediate time interval, which we will attempt to quantify in Sec.~\ref{sec:eigenvalues}, a different adjustment regime would be expected, which may or may not be self-similar. 

Next, using the Galerkin approach developed in Sec.~\ref{sec:ss_numerical}, we find numerical solutions of IBVP~\eqref{eq:IBVP_6} with $u^0(x) = \delta(x)$, for $\Bo\ne0$. Figure~\ref{fig:delta_function_IC}(c,d) presents the profiles of the rescaled numerical solutions for the representative values $\Bo = 1000$ and $3000$, respectively --- in terms of self-similar variables as $t^{1/6} u(\zeta)$ vs.\ $\zeta = x/t^{1/6}$ --- alongside the infinite-line self-similar profile $\mathcal{G}$ for $\Bo=0$. The selected plotting times are the same as those in Fig.~\ref{fig:delta_function_IC}(a,b). For both nonzero values of $\Bo$, we observe a reasonable collapse onto the self-similar profile at early times, similar to the case of $\Bo=0$ in Fig.~\ref{fig:delta_function_IC}(b). However, beyond $t \gtrsim 10^{-7}$, the rescaled profile for $\Bo = 1000$ in Fig.~\ref{fig:delta_function_IC}(c) ``drifts'' away from the self-similar profile $\mathcal{G}$, as it is no longer an exact solution to this problem. For the larger value of $\Bo=3000$ in Fig.~\ref{fig:delta_function_IC}(d), the rescaled numerical solutions further separate from the attractor.

\begin{figure}[ht]
	\centering
	\includegraphics[width=\textwidth]{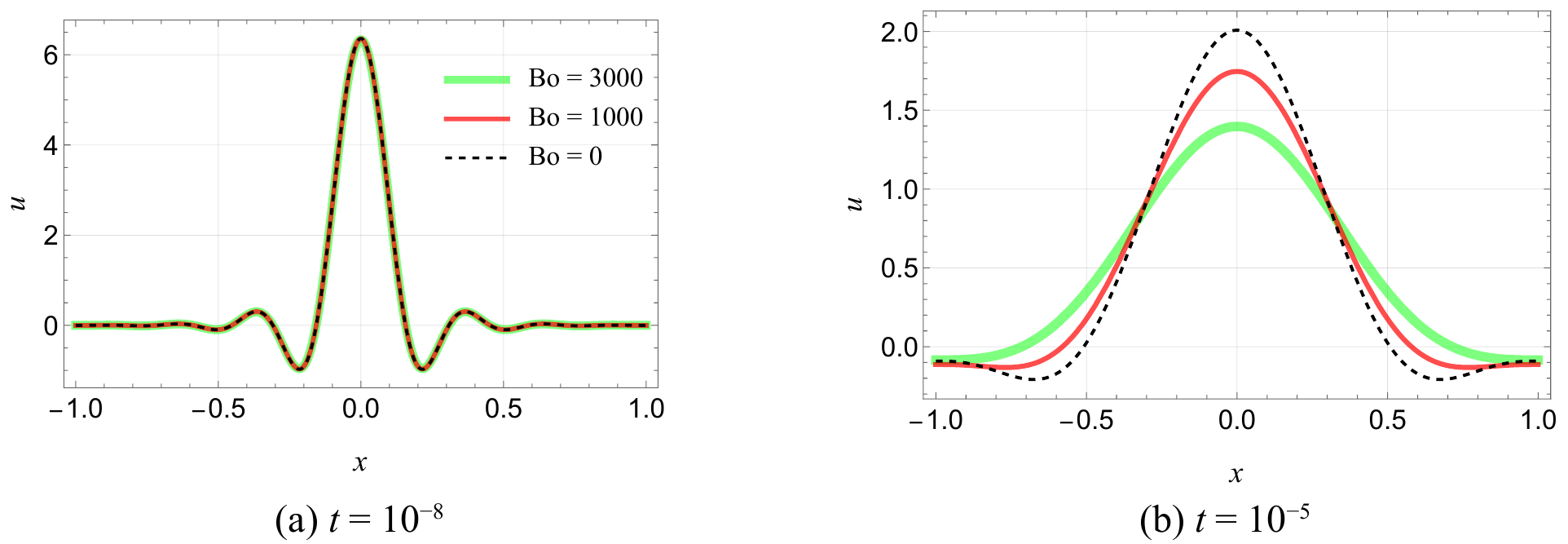}
	\caption{Snapshots of solutions $u(x,t)$ of IBVP~\eqref{eq:IBVP_6} with $u^0(x) = \delta(x)$ for different $\Bo$ and (a) $t=10^{-8}$ and (b) $t=10^{-5}$. The green (bright) and red (dark) solid lines are the solutions for $\Bo=3000$ and $\Bo=1000$, respectively, obtained numerically from the Galerkin expansion~\eqref{eq:u_x_expansion_formula}, whereas the black dashed curve is the case $\Bo=0$ obtained from the exact series solution~\eqref{eq:Green_2D_interval_convl} with coefficients~\eqref{eq:delta_exp}. For all presented cases, the spectral series was truncated at $M=50$.} 
	\label{fig:diff_G_delta_function_IC}
\end{figure}

Finally, to make the effect of $\Bo$ more clear, Fig.~\ref{fig:diff_G_delta_function_IC}(a,b) shows the unscaled numerical solution ($u$ vs.\ $x$) at two instants of time, early and late, respectively, for three different values of $\Bo$ each. It is now evident that the effect of $\Bo\ne0$ is weak at early times, as in Fig.~\ref{fig:diff_G_delta_function_IC}(a) in which all curves reasonably overlay the $\Bo=0$ solution. However, at late times, as in Fig.~\ref{fig:diff_G_delta_function_IC}(b), there is considerable spread in the solutions caused by $\Bo\ne0$. Further, while the profiles are indistinguishable in Fig.~\ref{fig:diff_G_delta_function_IC}(a) for the early time of $t=10^{-8}$, the profiles have separated at the later time of $t=10^{-5}$ in Fig.~\ref{fig:diff_G_delta_function_IC}(b), and the profiles for larger $\Bo$ have decayed more.

\subsubsection{Stepped initial condition}\label{sec:sss_results_step}

\begin{figure}[ht]
	\centering
	\includegraphics[width=\textwidth]{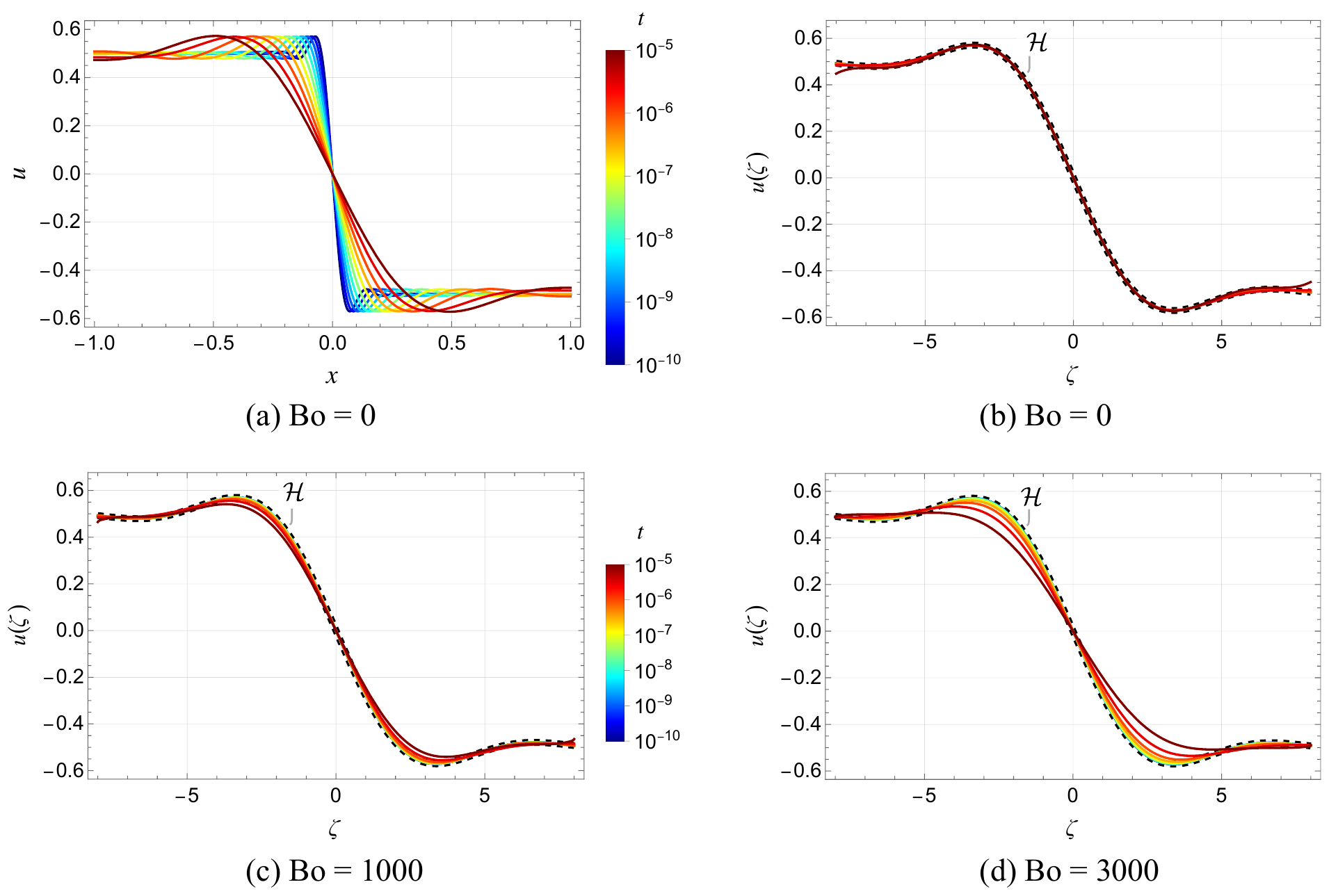}
	\caption{(a) Time-evolution of the solution $u(x,t)$, given by Eqs.~\eqref{eq:Green_2D_interval_convl} and \eqref{eq:odd_step_defn}, for $u^0(x) = H_\mathrm{odd}(x)$ with $\Bo=0$. (b) The profiles from (a) under the self-similar rescaling $u \mapsto u(\zeta)$, $\zeta = x/t^{1/6}$. The thick dashed black curve indicated with $\mathcal{H}$ is the self-similar attractor \eqref{eq:attractor_2D_line_odd_step} obtained using the convolution formula~\eqref{eq:Green_2D_line_convl} based on the inifinte-line Green function. (c,d) The equivalent of (b), but for the numerical solution $u(x,t)$ of IBVP~\eqref{eq:IBVP_6} with $u^0(x) = H_\mathrm{odd}(x)$ for $\Bo=1000$ and $3000$, respectively.}
	\label{fig:Odd_Step_IC}
\end{figure}

Now, we turn to the solution of IBVP~\eqref{eq:IBVP_6} with $u^0(x)=H_\mathrm{odd}(x)$ for $\Bo=0$. Figure~\ref{fig:Odd_Step_IC}(a,b) shows the profiles of the finite-interval series solution~\eqref{eq:Green_2D_interval_convl} with coefficients~\eqref{eq:odd_step_exp} for eleven different times $t\in[10^{-10},10^{-5}]$, as before.

\citet{Pedersen2021} solved the axisymmetric problem with the initial condition $H(1-r)$ in terms of the planar radial coordinate $r$. However, a feature of the axisymmetric problem is that the initial condition $H(1-r)$ is actually a localized cylinder, so it eventually also converges on the point-load Green's function (the self-similar profile from the previous subsection). At \emph{very} early times, however, \citet{Pedersen2021} found that their stepped initial condition exhibits different self-similar behavior, given by a suitable convolution of their Green's function and the initial condition. The solution profiles shown in Fig.~\ref{fig:Odd_Step_IC}(a) resemble those in the main part of Fig.~2 of \cite{Pedersen2021}, and the rescaled solutions (in terms of the self-similarity variable) shown in Fig~\ref{fig:Odd_Step_IC}(b) are also similar to those in the inset of Fig.~2 of \cite{Pedersen2021}. When making our observations, we must take into account that, here, the time range presented ($t\in[10^{-10}, 10^{-5}]$) includes intermediate time values at which the solution has decayed considerably away from the IC, whereas Fig.~2 of \cite{Pedersen2021} concentrates on earlier times ($t\in[10^{-13}, 10^{-9}]$ approximately) to observe this behavior. 

Indeed, on a finite interval, we do not have the problem that the stepped initial condition is actually a cylinder, and we easily observe that the series solution~\eqref{eq:Green_2D_interval_convl} with coefficients~\eqref{eq:odd_step_exp} converges to the self-similar profile $\mathcal{H}$ from Eq.~\eqref{eq:attractor_2D_line_odd_step}. As we consider a finite domain, the effects of the boundaries are eventually ``felt'', and the solution drifts away from the attractor at late times, especially for larger $|\zeta|$. Quantifying this transition is the subject of the next subsection.

\subsubsection{Leading eigenvalues and adjustment times}\label{sec:eigenvalues}

To understand how the solutions to our model problems ``relax'' towards the equilibrium state $u=const.$  for $\Bo\ge0$, it is instructive to study the eigenvalues of the dynamical systems~\eqref{eq:gen_Matr_DS_DEs}, which we observe are autonomous. Further, the matrices $\bm{\mathrm{A}^c}$ and $\bm{\mathrm{A}^s}$ are constant, real, and symmetric, meaning that they are diagonalizable with $M$ real eigenvalues, which we index by $j$. Specifically, let $r_j^c$ and $\bm{\chi^c}_j$ be the eigenvalues and associated eigenvectors of $\bm{\mathrm{A}^c}$, while $r_j^s$ and $\bm{\chi^s}_j$ are the eigenvalues and associated eigenvectors of  $\bm{\mathrm{A}^s}$. Then, the \emph{fundamental} solutions of Eqs.~\eqref{eq:gen_Matr_DS_DEs_even} and \eqref{eq:gen_Matr_DS_DEs_odd}, can be expressed respectively as 
\begin{subequations}
	\label{eq:soln_form_DS} 
	\begin{align}
		\bm{\phi^c}_j(t) &= \bm{\chi^c}_j \me^{r_j^c t}, \qquad (j=1,2,\ldots,M) \\
		\bm{\phi^s}_j(t) &= \bm{\chi^s}_j \me^{r_j^s t}, \qquad (j=1,2,\ldots,M) 
	\end{align}
\end{subequations}
and all solutions $\bm{\hat{u}^c}(t)$, $\bm{\hat{u}^s}(t)$ of the dynamical systems~\eqref{eq:gen_Matr_DS_DEs_even} and \eqref{eq:gen_Matr_DS_DEs_odd} are linear combinations of the $\bm{\phi^c}_j(t)$ and $\bm{\phi^s}_j(t)$. The symmetry of the matrices implies that this representation holds even if some eigenvalues were to have algebraic multiplicity greater than one \cite{BoyceDiPrima_10thEd}. From Eqs.~\eqref{eq:soln_form_DS}, it is evident that the eigenvalues $r_j^c$ and  $r_j^s$ determine the evolution of the solution $u$ of the original IBVP. The leading (largest) eigenvalues, $r_1^c$ and $r_1^s$, and the associated modes $\bm{\phi^c}_1$ and $\bm{\phi^s}_1$, are expected to dominate the solution behavior (i.e., set the main features of the shape of the profile $u$). If these modes decay (i.e., the leading eigenvalue is negative), the solution will swiftly progress toward the equilibrium $u=const.$ allowed by the chosen BCs.

\begin{figure}
	\centering
	\includegraphics[width=\textwidth]{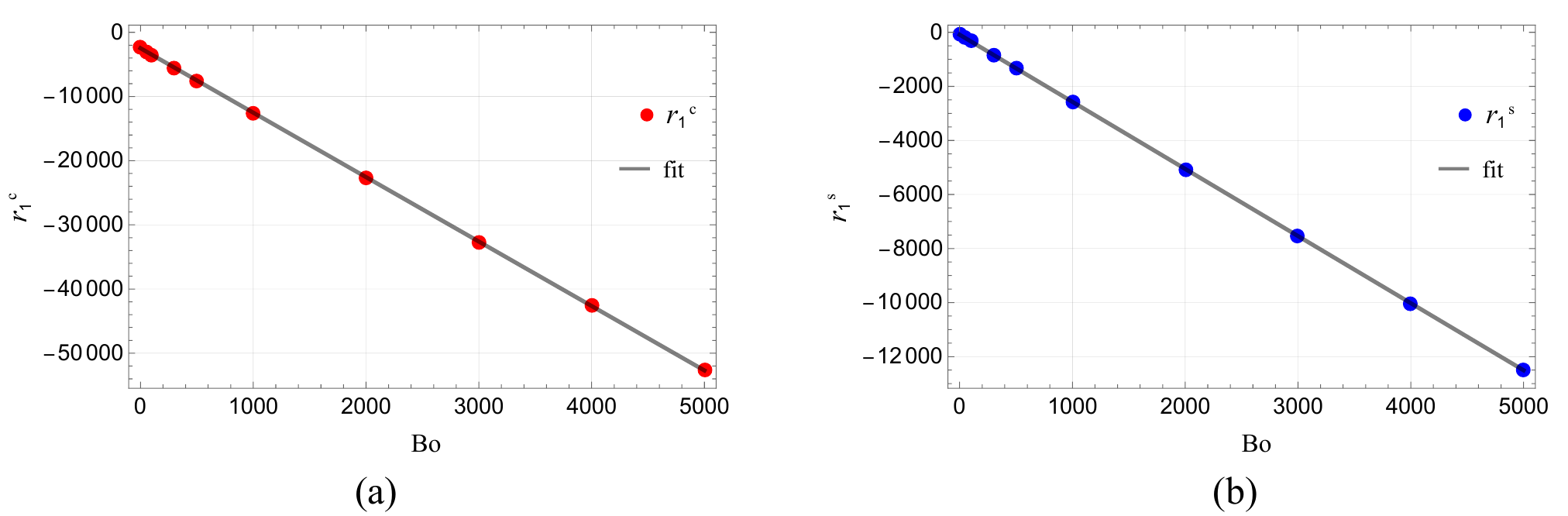}
	\caption{(a) Dominant eigenvalue $r_1^c$ (symbols) of the matrix $\bm{\mathrm{A}^c}$ as a function of $\Bo$ and the best-fit line $r^c_1(\Bo) = -(\lambda_{1}^c)^6 - 10.043915648\,\Bo$, where $\lambda_{1}^c = 3.66606496814$ \cite{NectarIvan2023}. (b) Dominant eigenvalue $r_1^s$ (symbols) of the matrix $\bm{\mathrm{A}^s}$ as a function of $\Bo$ and the best-fit line $r^s_1(\Bo) = -(\lambda_{1}^s)^6 - 2.483518333\,\Bo$, where $\lambda_{1}^s = 2.07175679767$ \cite{NectarIvan2023}. For both panels, the spectral series was truncated at $M=50$.}
	\label{fig:eigenvalues_DS}
\end{figure}

Therefore, it is interesting to investigate the effect of the model's parameters, specifically $\Bo$, on these eigenvalues. First, we focus on the even autonomous system formed by Eqs.~\eqref{eq:gen_Matr_DS_DEs_even} and \eqref{eq:gen_DS_Matr_defn_even}. We varied the value of $\Bo$. The eigenvalues of $\bm{\mathrm{A}^c}$ were computed using \textsc{Mathematica}. We considered a large range of $\Bo$ to observe the asymptotic scaling of eigenvalues with $\Bo$; at small $\Bo$ values, $\bm{\mathrm{A}^c}$ is dominated by the addends $-(\lambda_{m}^c)^6$ on its diagonal. We found that the eigenvalues were distinct for all the examined values of $\Bo$, but more importantly, \emph{negative}. The importance of this is twofold: Firstly, the solution of the dynamical system and, therefore, the IBVP decays with time, which is in agreement with the expected physical behavior. Secondly, from the numerical stability point of view, having all negative eigenvalues allows for the use of higher-order BDF schemes (such as BDF4 and BDF6), which are zero-stable for the entire negative real axis \cite{SuliMayers2006}.  Repeating the same analysis for the odd autonomous system formed by Eqs.~\eqref{eq:gen_Matr_DS_DEs_odd} and \eqref{eq:gen_DS_Matr_defn_odd}, we arrive at similar conclusions. Figure~\ref{fig:eigenvalues_DS} demonstrates the effect of $\Bo$ on the eigenvalues  (a) $r_1^c$ and (b) $r_1^s$, which correspond to the most dominant modes of the even and odd dynamical systems, respectively. Interestingly, it is observed that both $r_1^c$ and $r_1^s$ decrease linearly with $\Bo$ (best fits provided in the figure caption), which is different from the scaling of the leading eigenvalues with the Bond number obtained by \citet{Gabay2023} for a confined capillary thin film without interface elasticity.

\begin{figure}[ht]
	\centering
	\includegraphics[width=\textwidth]{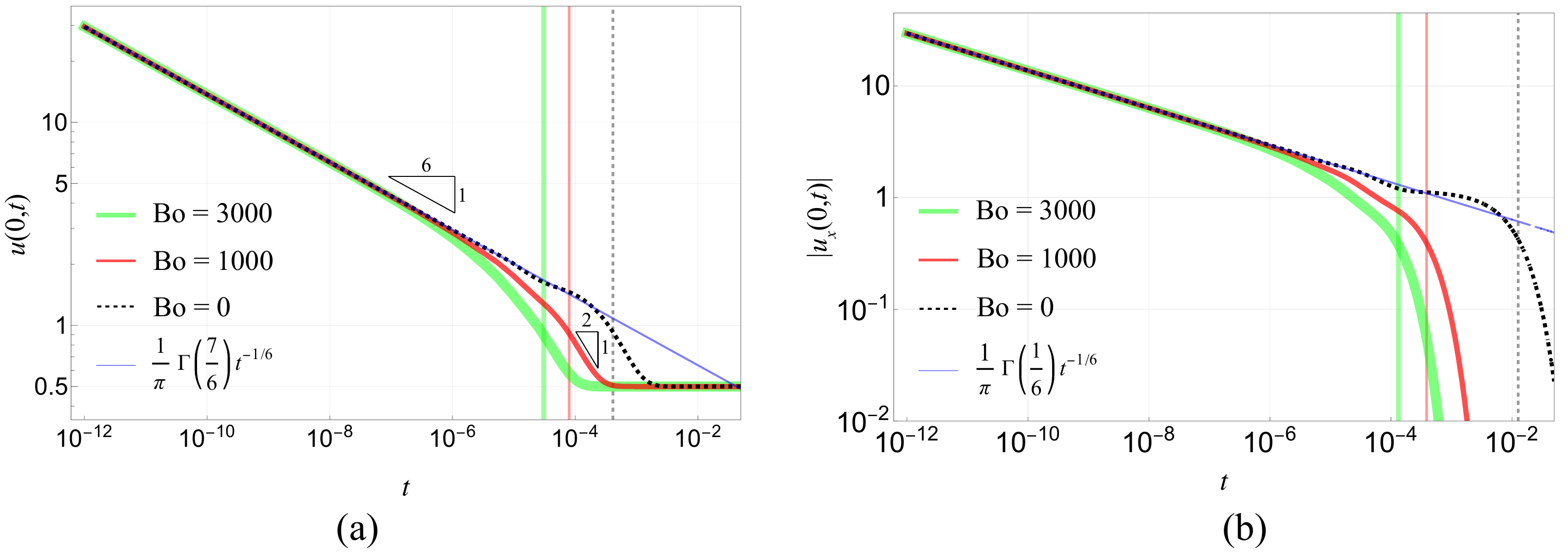}
	\caption{Estimating the time to diverge from the self-similar behavior. (a) Log-log plot of the maximum value $u(0,t)$ of the solution of IBVP~\eqref{eq:IBVP_6} with $u^0(x)=\delta(x)$. (b) Log-log plot of the maximum absolute value $|u_x(0,t)|$ of the slope of the solution of  IBVP~\eqref{eq:IBVP_6} with $u^0(x)=H_\mathrm{odd}(x)$. The thick light (green) and medium dark (red) solid curves are for $\Bo=3000$ and $\Bo=1000$, respectively, whereas the dashed black curves are for $\Bo=0$. The thin solid trendlines represent the expected self-similar scaling of the form $c\,t^{-1/6}$. The vertical lines of corresponding style and color (associated with the value of $\Bo$) are the estimated times, $t_e$, at which the solutions are expected to diverge from the self-similar behavior. For all presented cases, the spectral series was truncated at $M=50$.}
	\label{fig:decay_time}
\end{figure}

The decay of the solution towards the equilibrium is controlled by the leading eigenvalues of the operator. From the physical point of view, a characteristic equilibration time can be estimated, following, e.g., \cite{Liu2023}, as $t_e = -1/r_1^{c}$ or $-1/r_1^{s}$. To assess the physical meaning of $t_e$, in Fig.~\ref{fig:decay_time}(a), we show the maximum value $u(0,t)$ of the solution to IBVP~\eqref{eq:IBVP_6} with $u^0(x)=\delta(x)$, for different values of the Bond number $\Bo$, whereas Fig.~\ref{fig:decay_time}(b) depicts the maximum absolute value $|u_x(0,t)|$ of the slope of the solution to IBVP~\eqref{eq:IBVP_6} with $u^0(x)=H_\mathrm{odd}(x)$. Consistent with the self-similar behavior of these solutions, the trendline curves shown in Fig.~\ref{fig:decay_time}(a,b) are both of the form $c\,t^{-1/6}$, for some constant $c>0$. For (a), $c=\mathcal{G}(0)=\frac{1}{\pi}\Gamma(\frac{7}{6})$ via Eq.~\eqref{eq:attractor_delta}, whereas for (b), $c=|\mathcal{H}'(0)|=\frac{1}{\pi}\Gamma(\frac{1}{6})$ via Eq.~\eqref{eq:attractor_2D_line_odd_step}, where the prime denotes differentiation with respect to $\zeta$. As observed in Fig.~\ref{fig:decay_time}(a,b), both $u(0,t)$  and $|u_x(0,t)|$ decay as $t^{-1/6}$ up to approximately $t=10^{-6}$, for all three examined values of $\Bo$. 
		
Recall that, in Sec.~\ref{sec:sss_results_delta} (specifically, Fig.~\ref{fig:diff_G_delta_function_IC}) and Sec.~\ref{sec:sss_results_step}, we found that as time progressed, the solutions separated from the self-similar profile, with the solutions corresponding to larger $\Bo$ values separating further, and at earlier times. In this prior discussion, we limited our attention to times $t\le 10^{-5}$, and indeed, we see that all curves (for all $\Bo$ values considered) in Fig.~\ref{fig:decay_time} have left the self-similar trendline by about $t\simeq 10^{-5}$. Further, in the context of Fig.~\ref{fig:decay_time}, the time evolution towards and/or away from the attractor is determined by the convergence of the figure of merit onto and/or away from the trendline. Again, by our choice of initial conditions, it is evident that the solutions start \emph{on} the attractor, but eventually (as time increases), the figure of merit (solution or solution slope at the origin) deviates from the trendline. The characteristic time $t_e$ was calculated for each value of $\Bo$ considered and is shown as a vertical line in Fig.~\ref{fig:decay_time} and, according to this figure, reasonably predicts whether the solution has left the self-similar attractor, though it falls somewhere in-between the first deviation from the self-similar trendline and the flat portion of the curves beyond $t\gtrsim 10^{-4}$ when leveling on the finite interval occurs. Moreover, the deviation occurs earlier for larger values of $\Bo$, consistent with the scalings of $r_1^{c}$ and $r_1^{s}$ with $\Bo$ found above. 

Following up on the last observation, for $\Bo$ large enough, one might expect that the second-derivative term in Eq.~\eqref{eq:PDE_6} would dominate and a transition from the $u(x,t) \simeq \mathcal{G}(x/t^{1/6})/t^{1/6}$ self-similar asymptotics may occur, to a regime with a ``classical'' diffusion scaling, such that $u(x,t) \simeq  \mathcal{G}(x/t^{1/2})/t^{1/2}$. Following the approach of \citet{cs12}, we may seek a solution (of each form) by  substituting each ansatz into Eq.~\eqref{eq:PDE_6}:
\begin{subequations}\begin{alignat}{3}
& u(x,t) = t^{1/6} \mathcal{G}(x/t^{1/6}) &\qquad\Rightarrow\qquad& 6\mathcal{G}'''''' + \underbrace{6\,\Bo\,t^{2/3}\mathcal{G}''}_{\ll1\text{ for }t\ll1} + (\zeta \mathcal{G})' &= 0,\\
& u(x,t) = t^{1/2} \mathcal{G}(x/t^{1/2}) &\qquad\Rightarrow\qquad& \underbrace{\frac{2}{t^2}\mathcal{G}''''''}_{\ll1\text{ for }t\gg1} + 2 \,\Bo\, \mathcal{G}'' + (\zeta \mathcal{G})' &= 0,
\end{alignat}\label{eq:approx_ss}\end{subequations}
where $'=\rd/\rd \zeta$. Since $t$ cannot be eliminated from the above equations, a self-similar solution does not exist. However, in the spirit of \cite{cs12}, we may argue that in the first case, the term preventing self-similarity ($\sim t^{2/3}$) can be neglected for $t\ll1$. And vice versa, in the second case, the term preventing self-similarity ($\sim t^{-2}$) can be neglected for $t\gg1$. Thus, one can expect self-similar \emph{asymptotics} for \emph{early} and \emph{late} times (suitably defined). However, since the problem is dimensionless and all scales have been eliminated, there are no remaining free parameters to determine a cross-over time. As the triangle slopes in Fig.~\ref{fig:decay_time}(a) suggest, the time at which the solution transits from the $u(0,t) \sim t^{1/6}$ to $u(0,t) \sim t^{1/2}$ is quite close to $t_e = -1/r^c_1(\Bo)$. However, the figure suggests a similar transient even for $\Bo=0$ (albeit much shorter), and none of the values of $t_e$ are even above 1, which makes the balancing argument based on terms of Eq.~\eqref{eq:approx_ss} unconvincing. Further research is necessary to theoretically understand the potentially multiple self-similar behaviors and transitions of a finite thin film with elastic resistance in the presence of gravity.

\section{Conclusion}\label{sec:conclusion}

We studied a model sixth-order parabolic PDE on a finite interval, in the context of a thin fluid film between three impermeable walls topped by an interface with elastic bending resistance. Using a suitable eigenfunction expansion arising from an associated sixth-order, self-adjoint eigenvalue problem \cite{NectarIvan2023}, we explicitly constructed the Green's function for the model PDE. Using the Green's function, we obtained analytical solutions for two canonical problems: a point (localized) disturbance on the film and a stepped film. In both cases, we compared the intermediate asymptotics of the constructed series solutions to the infinite-line solution, which is \textit{de jure} self-similar. We found that the finite-interval solutions collapse onto the corresponding self-similar infinite-line attractors for intermediate times, despite the fact that the finite boundary conditions at the lateral walls drive the solution to the flat (equilibrium) state. Furthermore, we showed that the time it takes for the solutions to diverge from the attractor can be estimated in terms of the leading eigenvalue of the series expansion.

When gravity plays a role in the leveling of such an elastic-plated film, the governing PDE also presents a second-order derivative. In this case, rather than constructing another eigenvalue problem and rederiving the eigenfunction expansions, we used a Galerkin idea to obtain numerical solutions to the governing IBVP with gravity. By expanding this solution in terms of the same eigenfunctions as the no-gravity problem, we obtained dynamical systems for the time-dependent coefficients. These dynamical systems were integrated in time using \textsc{Mathematica}'s \texttt{NDSolve}. The ability to accurately expand singular and discontinuous initial conditions into a suitable basis using a relatively small number of modes was crucial to obtaining highly accurate solutions. A similar endeavor using finite difference/element methods would not be straightforward, requiring both a regularization of these singular or discontinuous initial conditions as well as fine spatial grids. On the other hand, Chebyshev pseudospectral methods would require the use of sophisticated preconditioners to overcome the ill-conditioning of high-order Chebyshev differentiation matrices \cite{Wang_SIAM,Hesthaven_SIAM}, as well as bordering six rows of the matrix to enforce the six boundary conditions of our problem~\cite{Boyd}. We also computed the eigenvalues of the resulting matrices in our Galerkin semi-discretization and found that the leading eigenvalue scales linearly with the Bond number, extending our estimate of the characteristic time to diverge from the self-similar attractor.

In summary, we demonstrated that the sixth-order self-adjoint eigenvalue problem introduced in \cite{NectarIvan2023} provides a suitable and effective framework for studying (\emph{both} exactly and approximately) a variety of thin film problems of elastic interfaces with bending resistance on a \emph{finite} domain, distinct from the bulk of prior literature focusing on the \emph{infinite}-line problem. In future work, it would be of interest to extend the proposed Galerkin spectral method based on the sixth-order eigenfunctions to the nonlinear version of the thin-film equation~\eqref{eq:elastic_film_0} (see, e.g., Eq.~(2.7) in \cite{NectarIvan2023}), and to explore the self-similar intermediate asymptotics of the nonlinear problem on a finite domain.

Our quantitative observations have applications to various studies interrogating the self-similar dynamics of thin films, for which an experiment can never occur on an infinite domain (e.g., \cite{Carlson2016,Kodio2017,Pedersen2019,Peng2020,Pedersen2021,Saeter2024}). Boundary effects are always present, and determining the time interval over which the self-similar behavior being observed is not affected by the confining boundaries is of interest. Similarly, the Green's function of the linearized problem is needed to solve forward and inverse problems of thin film actuation for microfluidics applications \cite{Tulchinsky2016,Rubin2017,Boyko2019,Gabay2023} such as patterning of surfaces. In addition, the time scales for the decay of perturbations and how they are affected by closed boundaries are crucial to understanding the shaping of thin films \cite{Gabay2023} for applications such as freeform diffractive optical elements \cite{Eshel2022}.

\paragraph{Acknowledgements} I.C.C.\ would like to acknowledge the hospitality of the University of Nicosia, Cyprus, where this work was initiated during his visits as a Fulbright U.S.\ Scholar (2022--2023). This paper is dedicated with respect and admiration to the memory of Prof.\ Stephen H.\ Davis, a pioneer in interfacial fluid mechanics \cite{Miksis2024}, including the present topic of thin films \cite{Oron1997}. I.C.C.\ recollects many insightful interactions with Prof.\ Davis, while I.C.C.\ was a graduate student (2007--2011) in the Department of Engineering Sciences \& Applied Mathematics at Northwestern University.

This work is based on our earlier paper \cite{NectarIvan2023}, appearing in the \textit{Proceedings of the 15\textsuperscript{th} Conference of the Euro-American Consortium for Promoting the Application of Mathematics in Technical and Natural Sciences (AMiTaNS 2023)}, from which material has been properly quoted throughout and referenced accordingly, and can be reproduced under the Creative Commons Attribution 3.0 license.

\paragraph{Author Contributions} \textbf{N.C.P.:} Conceptualization, Methodology, Software, Validation, Formal analysis, Investigation, Writing -- Original Draft, Writing -- review \& editing, Visualization. 
\textbf{I.C.C.:} Conceptualization, Methodology, Software, Validation, Formal analysis, Investigation, Writing -- Original Draft, Writing -- review \& editing, Funding acquisition.

\paragraph{Data availability} This study did not generate new data. The plots can be reproduced from the equations and descriptions in the main text.

\paragraph{Funding}
ICC's work on interfacial dynamics was partially supported by a Fulbright U.S.\ Scholar award from the U.S.\ Department of State and the U.S.\ National Science Foundation under grant CMMI-2029540.

\subsection*{Declarations}
\paragraph{Conflict of interest}
The authors declare no competing interests.

\addcontentsline{toc}{section}{References}
\bibliography{Sixth-Order-JEM.bib}

\clearpage 
\begin{appendices}
\section{The sixth-order orthonormal eigenfunctions}\label{sec:CON}

Here, for completeness and to make the present work self-contained, we reproduce the key details from \cite{NectarIvan2023} regarding the sixth-order eigenfunctions.

The self-adjoint sixth-order EVP, which represents a high-order Sturm--Liouville problem following the terminology of \citet{GM98,GM00}, is
\begin{subequations}
	\label{eq:6th_EVP}
	\begin{align}
		\label{eq:ODE_hg_6} 
		 -\frac{\rd^6 \psi}{\rd x^6} &= \lambda^6 \psi , 	\\
		\left.\frac{\rd \psi}{\rd x}\right|_{x=\pm1} = \left.\frac{\rd^2 \psi}{\rd x^2}\right|_{x=\pm1}&= \left.\frac{\rd^5 \psi}{\rd x^5}\right|_{x=\pm1}=0 ,
		\label{eq:6th_BCs}
	\end{align}
\end{subequations}	
has a countable set of solutions $\{ \psi_n(x) \}_{n=1,2,\hdots}$, which can be split into even (``cosine'') $\{\psi_n^c\}_{n=1,2,\hdots}$ and odd (``sine'') $\{\psi_n^s\}_{n=1,2,\hdots}$ ones, respectively:
\begin{subequations}\label{eq:efuncs_6_125}
	\begin{align}
		\psi_n^c(x) &= c^c_n \Bigg\{
		\tfrac{ 4\sin{\lambda_n^c}}{\cos{\lambda_n^c} - \cosh{\sqrt{3}\,\lambda_n^c}} \bigg[ \sin{\tfrac{\lambda_n^c}{2}} \cosh{\tfrac{\sqrt{3}\,\lambda_n^c}{2}} \cos{\tfrac{\lambda_n^c}{2}  x} \cosh{\tfrac{\sqrt{3}\,\lambda_n^c}{2} x} \nonumber\\ 
		&\qquad\quad -\cos{\tfrac{\lambda_n^c}{2}} \sinh{\tfrac{\sqrt{3}\,\lambda_n^c}{2}} \sin{\tfrac{\lambda_n^c}{2} x} \sinh{\tfrac{\sqrt{3}\, \lambda_n^c }{2}  x} \bigg] + \cos{\lambda_n^c x}\Bigg\}, \label{eq:efuncs_6_125_even} \\
		\psi_n^s(x) &= c^s_n \Bigg\{
		\tfrac{ 4\cos{\lambda_n^s}}{\cos{\lambda_n^s} +  \cosh{\sqrt{3}\,\lambda_n^s} } \bigg[ \sin{\tfrac{\lambda_n^s}{2}} \sinh{\tfrac{\sqrt{3}\,\lambda_n^s}{2}} \cos{\tfrac{\lambda_n^s}{2}  x} \sinh{ \tfrac{\sqrt{3}\,\lambda_n^s}{2}  x}    \nonumber\\ 
		&\qquad\quad -\cos{\tfrac{\lambda_n^s}{2}} \cosh{\tfrac{\sqrt{3}\,\lambda_n^s}{2}} \sin{\tfrac{\lambda_n^s}{2} x} \cosh{\tfrac{\sqrt{3}\,\lambda_n^s}{2}  x} \bigg] +  \sin{\lambda_n^s x}\Bigg\}, \label{eq:efuncs_6_125_odd} 
	\end{align}
where, for convenience, we have defined
\begin{align}
	c^c_n &:= 2 \sqrt{\tfrac{ \lambda^c_n}{d_n^c}} \left[\cos\lambda^c_n - \cosh \sqrt{3}\,\lambda^c_n \right], \label{eq:c_n^c_defn}\\
	c^s_n &:= 2 \sqrt{\tfrac{\lambda_n^s}{d^s_n}} \left[\cos\lambda_n^s + \cosh \sqrt{3}\,\lambda_n^s \right],\label{eq:c_n^s_defn}\\
	d^c_n &:= \sin{4\lambda^c_n} + 6\lambda^c_n \big( 2 - \cos{2\lambda^c_n} \big) + 2 \lambda^c_n \cosh{2\sqrt{3}\lambda^c_n} + 2 \sin{2\lambda^c_n} \cosh^2{\sqrt{3}\,\lambda^c_n} \nonumber\\
	&\phantom{=} + \big[ \sin{\lambda^c_n} -3 \sin{3\lambda^c_n} + 4\lambda^c_n \big( \cos{3\lambda^c_n} - 3\cos{\lambda^c_n} \big) \big] \cosh{\sqrt{3}\,\lambda^c_n}  \nonumber\\
	&\phantom{=} + 4\sqrt{3} \big( \cos{\lambda^c_n} - \cosh{\sqrt{3}\lambda^c_n} \big)\sin^2{\lambda^c_n} \sinh{\sqrt{3}\,\lambda^c_n},\label{eq:d_n^c_defn}\\
	d^s_n &:= 12 \lambda^s_n - 3\sin{2\lambda^s_n} - \sin{4\lambda^s_n} + 10 \lambda^s_n \cos{2 \lambda^s_n}\nonumber\\
	&\phantom{=} + \big( 2\lambda^s_n - \sin{2\lambda^s_n} \big) \cosh 2 \sqrt{3}\lambda^s_n - 4\sqrt{3}\cos^2\lambda_n^s \sinh(\sqrt{3} \lambda^s_n) \big( \cos\lambda^s_n + \cosh\sqrt{3} \lambda^s_n \big) \nonumber\\
	&\phantom{=} + 2 \big[ 4 \lambda^s_n (2 + \cos{2\lambda^s_n}) - 3 \sin{2 \lambda^s_n} \big] \cos{\lambda^s_n} \cosh{\sqrt{3} \lambda^s_n}.
	\label{eq:d_n^s_defn}
\end{align}%
\end{subequations}

The corresponding eigenvalues $\lambda^c$ and $\lambda^s$ satisfy the transcendental equations
\begin{subequations}\label{eq:evals_6_125}\begin{align}
		\cos{2 \lambda^c} + \sqrt{3}\sin\lambda^c \sinh{\sqrt{3}\,\lambda^c} - \cos\lambda^c \cosh{\sqrt{3}\,\lambda^c} &= 0,  \label{eq:evals_e_3_6}\\
		\sin{2 \lambda^s} + \sqrt{3}\cos\lambda^s \sinh{\sqrt{3}\,\lambda^s} + \sin\lambda^s \cosh{\sqrt{3}\,\lambda^s} &= 0,  \label{eq:evals_o_3_6}
\end{align}\end{subequations}
and form a discrete set. For $n=1,2,\hdots$, Eqs.~\eqref{eq:evals_6_125} have asymptotic solutions $\lambda_n^c \sim (n+1/6)\pi$ and $\lambda_n^s \sim (n-1/3)\pi$ as $\lambda\to\infty$, accurate to 12 and 11 digits, respectively, for $n=6$.

In addition, $\lambda_0^c=0$ is an eigenvalue of EVP~\eqref{eq:6th_EVP} with corresponding eigenfunction $\psi_0^c(x)=1$.
Importantly, any function $u(x)\in L^2[-1,1]$ can be expanded in this basis as in Eq.~\eqref{eq:u_x_expansion_formula} above (taking $M\to\infty$ therein) \cite{CodLev}.

The coefficient formulas for expanding the second derivatives of the eigenfunctions into a series of basis functions were derived in~\cite{NectarIvan2023} and are also quoted here for completeness. From Eq.~\eqref{eq:sec_der_even_def} (even eigenfunctions), given a fixed $n=1,2,3,\ldots$, and any $m=1,2,3,\ldots$, we have
\begin{equation}
	\label{eq:sec_deriv_even_form}
	\beta_{nm}^c =
	\begin{cases} 
		\frac{6 c_n^c c_m^c (\lambda _m^c)^3 (\lambda _n^c)^3}{(\lambda _m^c)^6 - (\lambda _n^c)^6}\left[ \frac{\lambda _m^c \sin{\lambda _n^c} \left(\cos{2\lambda _m^c}-\sqrt{3} \sin{\lambda _m^c} \sinh{\sqrt{3}\,\lambda _m^c} - \cos{\lambda _m^c} \cosh{\sqrt{3}\lambda _m^c}\right)}{\cos{\lambda _m^c}-\cosh{\sqrt{3} \lambda _m^c}}\right. &\\
		\qquad \qquad \qquad \qquad \left. +\frac{\lambda _n^c \sin{\lambda _m^c} \left(-\cos{2\lambda _n^c}+\sqrt{3} \sin{ \lambda _n^c} \sinh{\sqrt{3}\,\lambda _n^c}+\cos{\lambda _n^c} \cosh{\sqrt{3} \lambda _n^c}\right)}{\cos{\lambda _n^c}-\cosh{\sqrt{3} \lambda _n^c}}\right]\,\ \hfill \text{for} \  n\neq m,&\\[7mm]
		- \frac{\lambda _n^c {(c_n^c)}^2}{{{8 (\cos \lambda _n^c - \cosh \sqrt{3} \lambda _n^c)^2}}}    \bigg[\lambda _n^c \left(3 \cos{2\lambda _n^c} + \sinh^2{\sqrt{3}\lambda _n^c}+\cosh^2{\sqrt{3}\,\lambda _n^c} \right. &\\
		\left. \qquad \qquad \qquad  \qquad \qquad \qquad \quad  \,+ 4 \sqrt{3} \sin^3{\lambda _n^c} \sinh{\sqrt{3}\,\lambda _n^c} - 4 \cos^3{\lambda _n^c}\cosh{\sqrt{3}\,\lambda _n^c}\right)  & \\
		\qquad \qquad \ +\,2 \sin{\lambda _n^c} (\cos{\lambda _n^c}-\cosh{\sqrt{3} \lambda _n^c}) &\\
		\qquad \qquad \qquad \times \left(\sqrt{3} \sin{\lambda _n^c} \sinh{\sqrt{3}\lambda _n^c}+\cos{\lambda _n^c} \cosh{\sqrt{3} \lambda _n^c} -\cos{2 \lambda _n^c}\right)\bigg]\,\ \hfill \text{for} \  n= m,&
	\end{cases}
\end{equation}
where $c_n^c$ and $c_m^c$ are given by Eq.~\eqref{eq:c_n^c_defn}. Similarly, from Eq~\eqref{eq:sec_der_odd_def} (odd eigenfunctions), we have
\begin{equation}
	\label{eq:sec_deriv_odd_form}
	\beta_{nm}^s =
	\begin{cases} 
		\frac{6 c_n^s c_m^s (\lambda _m^s)^3 (\lambda _n^s)^3}{(\lambda _m^s)^6 - (\lambda _n^s)^6} \left[ - \frac{\lambda _m^s \cos{\lambda _n^s} \left(\sin{2\lambda _m^s}-\sqrt{3} \cos{\lambda _m^s} \sinh{\sqrt{3}\,\lambda _m^s} + \sin{\lambda _m^s} \cosh{\sqrt{3}\lambda _m^s}\right)}{\cos{\lambda _m^s}+\cosh{\sqrt{3} \lambda _m^s}}\right. &\\
		\qquad \qquad \qquad \qquad \left. +\frac{\lambda _n^s \cos{\lambda _m^s} \left(\sin{2\lambda _n^s}-\sqrt{3} \cos{ \lambda _n^s} \sinh{\sqrt{3}\,\lambda _n^s}+\cos{\lambda _n^s} \cosh{\sqrt{3} \lambda _n^c}\right)}{\cos{\lambda _n^c}-\cosh{\sqrt{3} \lambda _n^c}}\right]\,\hfill \text{for} \  n\neq m,&\\[7mm]
		\frac{\lambda _n^s {(c_n^s)}^2}{{{2 (\cos \lambda _n^s + \cosh \sqrt{3} \lambda _n^s )^2}}}     \bigg[\lambda _n^s \left(- \cos{2\lambda _n^s} + \sinh^2{\sqrt{3}\lambda _n^s}+\cosh^2{\sqrt{3}\,\lambda _n^s} \right. &\\
		\left. \qquad  \qquad \qquad \qquad \quad  \,- 4 \sqrt{3} \cos^2{\lambda _n^s}\sin{\lambda_n^s} \sinh{\sqrt{3}\,\lambda _n^s} + 4 \sin^2{\lambda _n^s}\cos{\lambda_n^s}\cosh{\sqrt{3}\,\lambda _n^s}\right)  & \\
		\qquad \qquad \ +\,2 \cos{\lambda _n^s} (\cos{\lambda _n^s}+\cosh{\sqrt{3} \lambda _n^s}) &\\
		\qquad \qquad \qquad \times \left(\sqrt{3} \cos{\lambda _n^s} \sinh{\sqrt{3}\lambda _n^s} - \sin{\lambda _n^s} \cosh{\sqrt{3} \lambda _n^s}-\sin{2 \lambda _n^s}\right)\bigg]\,\hfill \text{for} \  n= m,&
	\end{cases}
\end{equation}
where $c_n^s$ and $c_m^s$ are given by Eq.~\eqref{eq:c_n^s_defn}.

\clearpage

\section{Intermediate asymptotics of an IBVP with an even post with narrow support as IC}\label{sec:narrow_support}

To better understand how a solution approaches and then moves away from the self-similar attractor $\mathcal{G}$, as time progresses, we reexamine the solution to IBVP~\eqref{eq:IBVP_6} for $\Bo=0$ with a localized, but not Dirac $\delta$, initial condition. Specifically, we consider the following even post with narrow support
\begin{equation}
	u^0(x) = 	\begin{cases}
		 15, &\quad -\frac{\D 1}{\D 30}\leq x\leq \frac{\D 1}{\D 30}, \\[2mm]
		 0,  &\quad \text{otherwise} .
	\end{cases}
	\label{eq:IC_narrow_post}
\end{equation}
The coefficients $u_m^0$ used in the expansion of the initial condition~\eqref{eq:IC_narrow_post}, into the even eigenfunctions $\psi^c_m$ from Eq.~\eqref{eq:efuncs_6_125_even}, are
\begin{subequations}
	\label{IC_narrow_post_coeff_expr}
	\begin{align}
		u^0_0 &=1,\\
		u^0_m &=\frac{p^c_m}{\sqrt{q^c_m}}, \qquad\text{where}\\
		p^c_m &=  60 \Bigg\{\sin{\tfrac{\lambda^c_m}{30}} \left(\cos{\lambda^c_m} -\cosh{\sqrt{3}\lambda^c_m}\right)  \\
		&\qquad\qquad +\sin{\lambda^c_m} \Big[\sqrt{3} \left(\sin{\tfrac{29\lambda^c_m}{60}} \sinh{\tfrac{31 \lambda^c_m}{20 \sqrt{3}}} -\sin{\tfrac{31 \lambda^c_m}{60}} \sinh{\tfrac{29 \lambda^c_m}{20 \sqrt{3}}}\right) \nonumber\\
		&\qquad\qquad\qquad\qquad -\cos{\tfrac{31\lambda^c_m}{60}} \cosh{\tfrac{29\lambda^c_m}{20\sqrt{3}}} + \cos{\tfrac{29 \lambda^c_m}{60}} \cosh{\tfrac{31\lambda^c_m}{20 \sqrt{3}}}    \Big]      \Bigg\} \nonumber \\
		q^c_m &=\lambda^c_m\Bigg\{ \sin{4\lambda^c_m} - 6\lambda^c_m (\cos{2\lambda^c_m}-2) + 2 \lambda^c_m\cosh{2\sqrt{3}\lambda^c_m} + 2\sin{2\lambda^c_m}\cosh^2{\sqrt{3}\lambda^c_m}\\ 
		&\qquad  \qquad +\cosh{\sqrt{3}\lambda^c_m} \Big[\sin{\lambda^c_m} - 3\sin{3\lambda^c_m} + 4\lambda^c_m\left(\cos{3\lambda^c_m} - 3\cos{\lambda^c_m}\right)\Big] \nonumber \\ 
		&\qquad  \qquad   +4 \sqrt{3} \sin^2{\lambda^c_m} \sinh{\sqrt{3}\lambda^c_m} \left(\cos{\lambda^c_m} - \cosh{\sqrt{3}\lambda^c_m}\right)\Bigg\} ,  \qquad \qquad \qquad    m=1,2,\hdots \ \ . \nonumber
	\end{align}
\end{subequations}

\begin{figure}
	\centering
	\includegraphics[width=\textwidth]{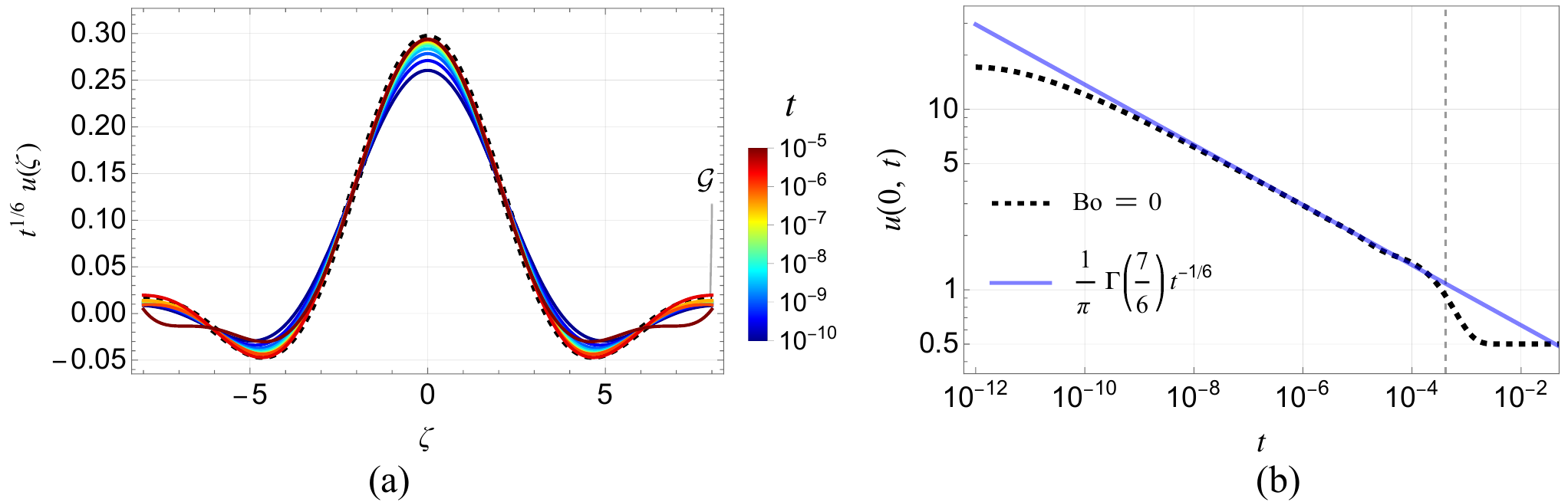}
	\caption{(a) The solution profiles of IBVP~\eqref{eq:IBVP_6} with IC~\eqref{eq:IC_narrow_post} for $\Bo=0$, under the self-similar rescaling, i.e., as $t^{1/6} u(\zeta)$ vs.\ $\zeta = x/t^{1/6}$. The selected times are the same as in Fig.~\ref{fig:delta_function_IC}. The thick dashed curve indicated with $\mathcal{G}$ is the self-similar attractor~\eqref{eq:attractor_delta} found via the infinite-line Green's function. (b)  Log-log plot of the maximum value $u(0,t)$ of the solution as a function of time $t$, for $\Bo=0$, as the dashed curve. The solid trendline represents the expected self-similar scaling of the form $c\,t^{-1/6}$. The black dashed vertical line is the estimated time, $t_e$, at which the solution is expected to diverge from the self-similar attractor. The spectral series was truncated at $M=50$.}
	\label{fig:narrow_post}
\end{figure}

Figure~\ref{fig:narrow_post}(a) presents snapshots of the numerical solutions (solid curves) for various times $t\in[10^{-10}, 10^{-5}]$, under the self-similar rescaling, i.e., as $t^{1/6} u(\zeta)$ vs.\ $\zeta = x/t^{1/6}$.  The self-similar attractor~\eqref{eq:attractor_delta} (thick dashed black curve) found via the infinite-line Green's function is shown on the same graph. From the figure, we observe that for early times, i.e., approximately $t\leq10^{-9}$, the solutions are at a significant distance from the self-similar profile $\mathcal{G}$, especially near $\zeta=0$. As time increases, we can see that the solution curves move towards the attractor, which demonstrates that when starting from non-ideal (i.e., neither a Dirac $\delta$ nor a unit step) initial conditions, the solutions are closest to the attractor for some \emph{intermediate} range of times, namely $10^{-8}\leq t \leq 5\times10^{-6}$ for this figure. For larger times, say $t=10^{-5}$, we can see that the overall shape of the solution has already changed compared to earlier times and is different from that of the attractor; notice the difference in concavity. 

Figure~\ref{fig:narrow_post}(b) depicts a log-log plot of the maximum value $u(0,t)$ of the solution as a function of time $t$ along with the self-similar trendline $\frac{1}{\pi}\Gamma(7/6)\,t^{-1/6}$. The estimated time $t_e = -1/r_1^c(0) = 1/\lambda_1^c$, at which the solution is expected to diverge from the self-similar attractor, is also marked on the graph. We observe that, in agreement with Fig.~\ref{fig:narrow_post}(a), the maximum value of the numerical solution is below the self-similar trendline for early times. Furthermore, it is \emph{attracted} to and is closest to the self-similar trendline for \emph{intermediate} times. Starting at $t_e \simeq -1/r_1^c(0) = 1/\lambda_1^c\simeq0.0004119$, the solution behavior moves away from the attractor as the effect of the boundaries becomes prominent, in agreement with what is observed in Fig.~\ref{fig:narrow_post}(a) for $t\geq10^{-5}$.

\end{appendices}

\end{document}